\documentclass[11pt]{article}
\pdfoutput=1
\usepackage[top=0.8in,bottom=1in,left=1in,right=1in]{geometry}
\usepackage{amsfonts,amsmath,amssymb,amsthm}
\usepackage{graphicx}
\usepackage{epstopdf}
\usepackage{algorithm, algorithmic}
\usepackage{bm}
\usepackage{caption}
\usepackage{subcaption}
\usepackage{cleveref}
\usepackage{authblk}

\crefformat{equation}{(#2#1#3)}

\ifpdf
  \DeclareGraphicsExtensions{.eps,.pdf,.png,.jpg}
\else
  \DeclareGraphicsExtensions{.eps}
\fi

\graphicspath{{./figures/}}

\title{A Fast Distributed Data-Assimilation Algorithm for Divergence-Free Advection\thanks{This research was partially
supported by US Department of Energy Contract
NDE-EE0006017. This report was prepared as an account of work sponsored by an agency of the United States government. Neither the United States government nor any agency thereof, nor any of their employees, makes any warranty, expressed or implied, or assumes any legal liability or responsibility for the accuracy, completeness, or usefulness of any information, apparatus, product, or process disclosed, or represented that its use would not infringe privately owned rights. Reference herein to any specific commercial product, process, or service by trade name, trademark, manufacturer, or otherwise does not necessarily constitute or imply its endorsement, recommendation, or favoring by the United States government or any agency thereof. The views and opinions of authors expressed herein do not necessarily state or reflect those of the United States government or any agency thereof.}}
\author[1]{Tigran T. Tchrakian\thanks{tchrakit@gmail.com}}
\author[1]{Sergiy Zhuk\thanks{sergiy.zhuk@ie.ibm.com}}
\affil[1]{IBM Research -- Ireland}
\date{}                     
\usepackage{amsopn}
\DeclareMathOperator{\diag}{diag}
\newtheorem{proposition}{Proposition}

\begin{document}

\maketitle

\begin{abstract}
In this paper, we introduce a new, fast data assimilation algorithm for a 2D linear advection equation with
divergence-free coefficients. We first apply the nodal discontinuous Galerkin
(DG) method to discretize the advection equation, and then employ a set of
interconnected minimax state estimators (filters) which run in parallel on spatial elements possessing
observations. The filters are interconnected by means of numerical Lax-Friedrichs fluxes. Each filter is
discretised in time by a symplectic Mobius time integrator which preserves all quadratic
invariants of the estimation error dynamics. The cost of the proposed algorithm scales linearly with the
number of elements. Examples are presented using both synthetic and real data. In the latter case, satellite images are
assimilated into a 2D model representing the motion of clouds across the surface of the Earth.
\end{abstract}

\section{Introduction}

Data Assimilation (DA) is an important component in many modern industrial cyber-physical systems. It improves the accuracy of forecasts provided by physical models by optimally combining their states -- \emph{a priori} knowledge encoded in equations of mathematical physics -- with \emph{a posteriori} information in the form of sensor data, and evaluates forecast reliability by taking into account uncertainty, i.e., model error or measurement noise. Mathematically, DA relies upon optimal control methods (deterministic state estimators) or applied probability (stochastic filtering). In the probabilistic framework, the optimal solution of the state estimation problem is given by the Kushner-Stratonovich equation, which describes the dynamics of the conditional probability density of the state of a Markov diffusion process given incomplete and noisy observations~\cite{GihmanSkorokhod1997}. In contrast, deterministic state estimators assume that errors have bounded energy and belong to a given bounding set. The state estimate is then defined as the minimax center of the reachability set, a set of all states of the physical model that are ``reachable'' from the given set of initial conditions and model errors, and are compatible with observations. The dynamics of the minimax center is described by a minimax filter~\cite{Kurzhanski1997}. In the case of linear models, the KS equation can be reduced to the Kalman-Bucy filter equations, which coincide with those of the minimax filter~\cite{KrenerMinimax}. We refer the reader to~\cite{ReichDA2015,StuartDA2015} for further discussion on modern data assimilation and state estimation.

A major issue of DA approaches is the lack of scalability: indeed, computing optimal state estimates for PDEs is often infeasible even in two spatial dimensions due to the ``curse of dimensionality''. To illustrate this issue consider a linear system
\begin{equation*}
  \dfrac{dc}{dt} = A c\,,\quad  y(t) = H c(t) +\eta(t)\,, \quad c(0) = c_0\,,
\end{equation*}
where $A$ is a differential operator representing the physical model, $y$ is the sensor data corrupted by some noise $\eta$, $H$ is the mathematical representation of the sensor relating states $c(t)$ to ${y}(t)$, and ${c}_0$ is an uncertain initial condition such that, for given symmetric positive definite operators $Q_0$ and $R$ it holds that: \[
(Q_0 {c}_0, {c}_0) \le 1\,, \quad \int_0^T E (R{\eta},{\eta}) dt\le 1\,,
\]
In this case, the optimal state estimator, $\hat {{c}}$, is given by the minimax filter:
\begin{align}
\dfrac{d\hat {{c}} }{dt} &= A \hat {{c}} + P H^\top R ({y} - H\hat{{c}})\,, & \hat {{c}}(0) = 0\,,\label{eq:flt}\\
\dfrac{dP}{dt} &= AP + PA^\star - P H^\star R H P\,, &P(0) = Q_0^{-1}\label{eq:ric}\,,
\end{align}
which is, in fact, equivalent to the estimate of the Kalman filter~\cite{KrenerMinimax}. As noted above, accurately approximating $\hat{{c}}$ and the gain (or state error covariance operator) $P$ in real time is not feasible even if $A$ is a linear advection-diffusion operator in two spatial dimensions: a very modest approximation of $A$ by $100$ or so basis functions in each spatial dimension will result in a $10000\times10000$ stiffness matrix so that, generally speaking, the problem of finding an accurate and fast approximation of $P$ solving the matrix differential Riccati equation~\cref{eq:ric} becomes intractable.

\subsection{Contributions}\label{sec:contributions}
In this paper, we propose an efficient and scalable data assimilation algorithm for linear advection equations with divergence-free coefficients in two spatial dimensions. The algorithm is distributed, i.e.,
it uses a network of local filters which process localised observations, and then exchange portions of the information with the neighbouring filters to reconstruct the ``global'' state. Mathematically, the proposed algorithm relies upon the nodal Discontinuous Galerkin (DG) discretization of the advection equation, and minimax state estimation framework. In what follows, we will briefly discuss the distributed filtering and its intrinsic relation to the DG discretization. This relationship is the corner-stone of the proposed DA algorithm.

\begin{figure}[bh]
\centering
                \includegraphics[width=\textwidth]{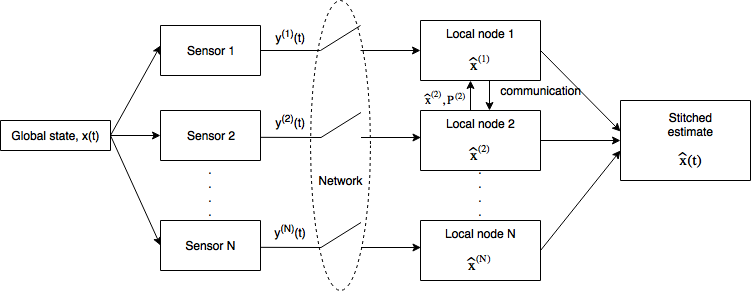}
                \caption{Distributed filtering strategy }
                \label{fig:dfilt}
\end{figure}

Distributed filtering/state estimation is popular in control engineering, specifically in distributed sensor networks where one of the typical problems is how to construct an estimate of the entire state (global estimate) of a dynamic spatio-temporal process from spatially distributed nodes of sensors in a decentralized way. This means that there is no ``central unit'' that runs a global model of the process and gathers all of the measurements from the nodes in order to construct a global state estimate by an appropriate technique, e.g.  Kalman filter. Instead, each node has a local model of the process, i.e. a model describing the dynamics of the restriction of the state of the global process onto the spatial region of the node, and local measurements of this restriction. In addition, each node may receive information (e.g., state estimates) from the neighbouring nodes. The neighbours are determined by the network topology. The goal of the distributed state estimation is to construct a local estimate at each node, i.e. the estimates of the restriction of the entire state to the node, by using all of the data available at the node. The local estimates must be constructed so that when they are stitched together, the resulting estimate of the entire state is as close as possible to the optimal global state estimate computed at the ``central unit'' as suggested above. The described distributed filtering strategy is visualised in~\cref{fig:dfilt}.\\
Note that the local nodes use local models of the process, and so these models must communicate with neighbours to maintain the ``global picture''. Thus, the key problem of distributed filtering is how to set up communications with the neighbours in order to have the ``stitched state estimate'' be as close as possible to the global state estimate. The key motivation for developing the distributed filtering is two-fold: on one hand, the global state estimation can be quite expensive computationally for reasons outlined above, and the nodes may not be in possession of sufficient computational power, and, on the other, communicating large amounts of information required to perform the global filtering may be very expensive and even infeasible especially if the communications are over a wireless network. We refer the reader to~\cite{olfati2005consensus} for further details on distributed state estimation strategies in the context of sensor networks.

Now, to relate the distributed filtering and the DG discretization we assume that the dynamic process of interest is modelled by a PDE, e.g. a linear advection equation, in a large spatial domain $\Omega$, and that the sensors are located only at some sub-domains of $\Omega$. 
It turns out that the DG discretization framework provides a natural means for implementation of the distributed filtering strategy, thus mitigating the computational cost associated with high dimensionality. Indeed, the local models of the spatio-temporal process are generated as follows: the DG method splits the computational domain $\Omega$ into a number of non-overlapping elements, $D^k$, and in each $D^k$, the original PDE is substituted by a local ODE. The latter has a local system matrix, $A^k$ which represents the advection operator $A$ in a local polynomial basis of $D^k$. All of the resulting local ODEs are interconnected by means of numerical flux functions which contribute to both local system matrices $A^k$ and source terms $\bm{b}^k$ to ``maintain the global picture''. Moreover, to compute the solution of the local ODE on the element $D^k$ at time level $t_s$ one needs only (i) the solution on $D^k$ in the past (e.g. at the previous time level, $t_{s-1}$), (ii) portions of solutions from adjacent elements (neighbours) in the past, and (iii) the advection field at time $t_s$ and possibly in the past. The data (i)-(iii) are used to update the matrices $A^k$ and source vectors $\bm{b^k}$, and (ii), in fact, implements the communication between the nodes by utilizing numerical fluxes. After the update, all of the local solutions are computed for time level $t_s$ independently. One of the key advantages of this discretization is that, by taking a large enough number of small elements $D^k$ with low order polynomial basis functions, one can approximate the solution of the PDE up to high precision, and, at the same time, keep the local system matrices small. This is in contrast to spectral/finite element discretization methods where one forms a global ``stiffness matrix'' to maintain the continuity of the solution (see~\cref{sec:state-art}).

As noted, the DG discretization framework naturally provides the ``localisation'' and ``communication'' components for the implementation of the distributed filtering strategy. Based on this, we propose a new distributed filtering  algorithm for advection equations by combining the DG discretization and minimax state estimation. Namely, for each element $D^k$, the local minimax state estimator (filter) is applied to estimate the state of the system of local ODEs by using the observations localised in $D^k$. The local filter relies upon a local Riccati operator $P^k$ which is constructed for every element $D^k$ by solving the standard Riccati equation given the local system matrix $A^k$, and the local (for $D^k$) uncertainty description (e.g. initial condition error, model error and observation noise). If $D^k$ has no observations, the local observation matrix, $H^k$ is set to zero, i.e. no data available on $D^k$. As a result, the innovation term equals zero, and the estimator coincides with the state of the local ODE. Note that in this case, $P^k$ has no impact on the dynamics of the estimator. 
The filters, located at adjacent elements, ``run'' independently and ``communicate'' with one another through numerical fluxes~(see \cref{sec:localMF}). In fact, the numerical flux is a function of the local state, the adjacent state and advection field restricted to the boundary between the elements. At every time step $t_s$, it updates the local system matrix $A^k$, which ``absorbs'' everything related to the solution on $D^k$, and the source term $\bm{b^k}$ which encapsulates all of the information about the state on adjacent elements. This provides implicit communications between local Riccati equations: indeed, solutions from adjacent elements ``change'' the solution on $D^k$ through the term $\bm{b^k}$ at time $t_s$, and this change is then reflected in the local system matrix $A^k$ after the numerical flux has been recomputed at the next time step $t_{s+1}$. In this way, implicitly, local Riccati matrices exchange information with each other.

The proposed distributed filtering algorithm has two notable advantages. One is its scalability, i.e. the potential for speed-up due to parallelisation of the filter over elements. Indeed, as noted above, the local solutions on each element $D^k$ can be obtained in parallel. As a result, the cost of the proposed algorithm scales linearly with the number of elements (see~\cref{sec:highResTest}). The other advantage is the possibility to preserve all of the quadratic invariants associated with local Riccati operators $P^k$. Indeed, the size of the system matrices is quite small, and this allows us to use a symplectic Mobius time integrator~\cite{ZhukCDC14} to discretize the corresponding matrix Riccati differential equation. This type of time discretization is implicit and preserves symmetry and positivity of $P^k$. In addition, all quadratic invariants of the estimation error dynamics are preserved too. The importance of this for practical applications (like the one in~\cref{sec:COD}) is that the simulation results are trustworthy and represent indeed what has been predicted by the theory for the continuous case.

Finally, we note that the localisation is an approximation, and as such it may lead to errors in approximating the solution of the Riccati equation, as the global gain, $P$ is informed by richer information than is $P^k$ in the distributed case. The localisation error is quantified by comparing the local estimates and operators $P^k$ versus \emph{the reference estimate} represented by optimal centralised minimax state estimate obtained by forming a ``global'' system matrix and solving the corresponding ``global'' Riccati equation. The resulting comparison shows that the distributed filter provides estimates that are quite close to the reference estimate (see~\cref{sec:globalFilter}).

\subsection{State-of-the-art}\label{sec:state-art}
In this paper the ``discretize-then-optimize'' strategy is used instead of discretizing an infinite dimensional filter/Riccati equation~\cite{Bensoussan} directly. The main reason for this is that the optimal infinite dimensional minimax/Kalman filters are available in the literature for linear parabolic equations including advection-diffusion equations, and symmetric systems of linear hyperbolic equations in $\mathbb{R}^n$. Neither suits our application, namely data assimilation of cloud optical depth images into an advection-dominated flow given by a scalar linear hyperbolic equations in a bounded domain with boundary conditions on the inflow zone. For this reason we adopt a hybrid strategy: discretize in space by using DG method, and perform the state estimation for the resulting continous time system of ODEs by employing the minimax filter.\\
Our choice of the DG method is motivated not only by its elemental nature which provides a natural means of implementation of the distributed filtering strategy, but also because of the stabilizing effect of the flux transfer between the elements, which makes it more robust in advection dominated settings~\cite{QiuShu, Cockburn2001, cockburn1989tvb}. This is a key difference between DG and the Spectral-Element Method (SEM) family of methods. Although the latter also use an elemental discretisation, their handling of element interfaces requires the field representing the advected quantity to be smooth (at least of $H^1$-class) across interfaces which can result in instabilities for advection-dominated flows~\cite{Malm2013} even with smooth advection coefficients. In order to overcome this issue, artificial diffusion can be added to the model~\cite{Malm2013}, but in many cases, the amount of diffusion required may be un-physical. We stress that in this paper we will restrict our attention to 2D advection by divergence-free (volume preserving) velocity fields, and thus will not need to apply any post-processing to our estimates. However, in general, for nonlinear problems and linear advection problems where the velocity field is non-smooth, a postprocessing in the form of a flux limiter can be employed for DG~\cite{QiuShu} to get convergent numerical approximations.\\
The reason for using the minimax filter is in that it assumes merely bounded model errors in contrast to the standard statistical assumptions of the Kalman-Bucy filter. Since we have to take the discretization error into account, and there is no statistical information about its distribution, the worst-case estimation error of the minimax filter appears to be more appropriate here. Note that, in the linear case, the state estimates provided by the Kalman-Bucy and minimax filters coincide~\cite{KrenerMinimax} -- the difference lies only in the interpretation, i.e. conditional mean versus minimax center of the reachability set.

The problem of finding approximations of the solution $P$ of the matrix differential Riccati equation has been studied by many authors. The most popular statistical method is the Ensemble Kalman Filter (EnKF)~\cite{EvensenOceanDyn2004}. Using that approach the computational bottleneck in computing $P$ is overcome by generating an ensemble of trajectories, e.g. ensemble of grid-functions each of which approximates the solution over the entire domain $\Omega$, and by computing the ensemble variance to approximate the state error covariance matrix $P$. The latter is then used to compute a state estimate in the same way as in the Kalman filter. We stress that if the size of the discretized state vector (total number of grid points) is quite large then, to achieve convergence of the empirical moments to the exact ones the number of ensemble members must be very large that is hardly feasible even for the advection-diffusion equation in two dimensions. To overcome this issue, the localised EnKF, a combination of the square-root Kalman Filter with EnKF, was proposed~\cite{ott2004local}. Similarly to what has been suggested above, the local EnKF splits the domain into regions $E^k$, e.g. rectangles, then it makes a localisation, namely the local ensemble vectors are created by restricting global ensemble vectors (grid-functions representing the solution over the entire domain at time $t_s$) to the regions $E^k$, and the local (empirical)  covariance matrices $P^k$ are computed on each $E^k$ from local ensembles. Given the specific structure of $P^k$ one can relatively easily compute its eigenvectors and eigenvalues and then perform the filtering step in the space spanned by the leading eigenvectors corresponding to the largest eigenvalues. This provides an update for the local ensembles. Finally, to go to the next time step $t_{s+1}$, the updated global ensemble vectors can be obtained from the updated local ensemble vectors by a weighted average. This algorithm does scale very well as all the local ensembles are completely independent. i.e. the local filtering steps in each region are not influenced by their neighbours. In contrast, the local state estimates generated by our algorithm do depend on the adjacent local estimates by means of the boundary integral interconnection mechanism of the DG method (see \cref{sec:localMF}) and so do the local gains. This mechanism exchanges the information ``in the direction'' of the advection coefficients and hence it maintains numerical stability without the need for artificial smoothing. In particular, it allows discontinuous observations to be dealt with (e.g. when every other element has observations, giving rise to a ``chequered'' pattern, see~\cref{sec:sparseObsSynth}) without the need to average the local estimates.

Finally, for the sake of completeness we mention here a family of Particle Filters (PF). Theoretically, one can use PFs to approximate the solution of the Zakai equation, which describes the dynamics of the (un-normalized) density of the conditional distribution of the (hidden) Markov diffusion given observations~\cite{crisan1997nonlinear}, or to construct a Gibbs sampler with a target Feynman-Kac distribution~\cite{del2014particle}. PF is an active research area and, at the moment, its application to state estimation problems for PDEs is quite limited~\cite{snyder2008obstacles}. We refer the reader to~\cite{rebeschini2015can} where a theoretical framework for localization of PFs is provided, and to~\cite{poterjoy2016localized} where an attempt to design an algorithm for localization of PF has been presented together with application to a small-to-medium sized state estimation problem for L96 model.

\subsection{Experimental assessment} The practical motivation of this work stems from the following problem. For accurate short-term solar energy forecasting (minutes- to hours-ahead) of large geographical areas (continental scale), models using geostationary operational environmental satellite (GOES) imagery are more effective than numerical weather prediction models as the latter typically take too long (several hours) to ramp up. The GOES satellite imagery only provides information on the current distribution of clouds, including cloud optical depth (COD) and top/bottom altitude, which can be converted to the current solar irradiance at the surface of the earth using radiative transfer modeling. An accurate cloud advection model is thus required to forecast the future cloud distribution and solar irradiance. We refer the reader to~\cite{SZTTAA_CDC17,zhuk2017computer} for further details on this. In~\cref{sec:experiments} we apply the proposed method for cloud advection over the domain with constant inflow and free-exit boundary condition. The domain spans 16 degrees of longitude and 12 degrees of latitude, and covers approximately $1.68e6$ square kilometres.

In addition, we consider a synthetic scenario in which a non-stationary, divergence-free velocity field is advecting a smooth quantity over a domain with non-stationary boundary conditions. 
We use a relatively low-resolution grid here, as the purpose of this scenario is comparison between ``local'' and ``global'' filters, with the latter being too computationally expensive to run at high resolution. In the case where observations are incomplete, elements with data neighbouring elements without data may give rise to sharp discontinuities, which will appear in the global system matrix. These discontinuities can cause numerical instabilities which can be avoided either by postprocessing, or by ``placing the discontinuities'' into the source term as in the case of the distributed filter. We implement the latter approach since the discontinuities are artificial, i.e. they have nothing to do with the physics of the advection process, and so ``placing them'' into the source terms is perfectly justified as they arise from the boundary terms imposing weak boundary conditions between the elements (see~\cref{sec:globalFilter}). 

Finally, we conduct a simple scalability study to demonstrate that the computational cost of the method scales linearly when the number of elements grows (up to $250\times250$) and the polynomial order is fixed ($N=3$) resulting in the discrete state vector of dimension $10^6$ (see~\cref{sec:highResTest}).

The paper is organized as follows. The problem statement is detailed in \cref{sec:probStat}, where we describe the advection model and the objectives of the filtering problem. Our main results are in \cref{sec:main} where we introduce and describe our filtering algorithm. Experimental results including comparisons of distributed and global filters for both synthetic and real data are  described in \cref{sec:experiments}, and the conclusions follow in \cref{sec:conclusions}. 

\paragraph{Notation}
$\mathbb{R}^n$ denotes $n$-dimensional vector space of real vectors $\bm{x}=(x_1\dots x_n)^\top$ with standard canonical basis and inner product $\bm{x}\cdot\bm{y}=\sum_{j=1}^n x_iy_i$. Let $\Gamma$ denote the boundary of the computational domain $\Omega$, $\Gamma^{in}$ -- the inflow zone of $\Gamma$, i.e. $\Gamma^{in}:=\{\bm{x}\in\Gamma: \bm{u}(x)\cdot\bm{n}(x)<0\}$, $\bm{n}$ is the unit outward vector which is normal to $\Gamma$ at $\bm{x}$, and $\bm{u}(\bm{x}) = (u(\bm{x}), v(\bm{x}))^\top$ is a given incompressible velocity field, i.e. $\nabla \cdot \bm{u} =0$. We also set $\Omega_T:=(0,T)\times\Omega$ and $\Gamma^{in}_T(0,T)\times\Gamma^{in}$ for a real number $T<+\infty$. Finally, let $L^2(\Omega)$ denote the space of square-integrable measurable functions over $\Omega$, and let $H^1(\Omega):=\{f\in L^2(\Omega):\partial_x f,\partial_y f \in L^2(\Omega)\}$.

\section{Problem Statement}
\label{sec:probStat}
The problem we aim to solve is that of state estimation or filtering for a linear advection equation in two spatial dimensions subject to uncertain but bounded error in the initial and boundary conditions, and uncertain forcing. In what follows we introduce the advection equation with uncertain parameters, provide the uncertainty description and formulate the problem statement.

\paragraph{State equation} 
Assume that the computational domain $\Omega \in \mathbb{R}^2$ is a convex bounded domain, and let the function $(\bm{x},t)\mapsto c(\bm{x},t)$ represent a quantity (e.g. a concentration of a material) that is advected by an incompressible velocity field\footnote{To simplify the presentation we assume that $\bm{u}$ is independent of time so that the inflow part of the boundary $\Gamma^{in}$ does not change over time. However, in some numerical examples, the proposed method is applied with time-dependent $\bm{u}$. } $\bm{u}$ according to the following linear advection equation:
\begin{equation}
\label{eq:advection}
\begin{split}
&\partial_t c(\bm{x},t)+ \bm{u}(\bm{x}) \cdot \nabla c(\bm{x},t) = g(\bm{x},t)e(\bm{x},t)\text{ in } \Omega_T\,,\\
&c(\bm{x} ,0) = c_0(\bm{x}) + g_0(\bm{x})e_0(\bm{x})\,, \text{ on } \Omega\,, \\
&c(\bm{x},t) = c_\partial(\bm{x},t) + g_\partial(\bm{x},t)e_\partial(\bm{x},t) \,, \text{ on } \Gamma^{in}_T\,,
\end{split}
\end{equation}
subject to an uncertain forcing $e$, boundary condition $c_\partial$ and initial condition $c_0$ with uncertain additive errors $e_\partial$ and $e_0$ such that:
\begin{equation}\label{eq:model_error}
|q_\partial(\bm{x},t)e^2_\partial(\bm{x},t)|\le 1\,,\quad |q_0(\bm{x}) e_0^2(\bm{x})|\le1\,,\quad |q(\bm{x},t)e(\bm{x},t)| \le 1\,,
\end{equation}
provided $q_0$, $q_\partial$ and $q$ are given smooth weighting functions such that $q_0,q_\partial, q>q^\star>0$ for all $(\bm{x},t)\in\Omega_T$, and the inequalities are understood for every $\bm{x}\in\Omega$ outside perhaps a set of measure zero. Note that $q_0$ and $q$ quantify the level of confidence in the initial condition/boundary conditions and state equation: namely, $q_0$ may specify ``zones'' of $\Omega$ where the knowledge of the initial condition $c_0$ is more precise or less so, and $q$ defines zones of $\Omega$ where~\cref{eq:advection} holds almost exactly ($|e|\approx 0$ in that zone) or only up to a significant error ($|e|>0$) and these zones may vary over time. Statistically, this corresponds to the maximal entropy assumption, i.e., any $(e_0,e,e_\partial)$ verifying~\cref{eq:model_error} have equal probability of appearing. Functions $g$, $g_0$ and $g_\partial$ are assumed to be given and allow one to either localize errors in space/time, or switch them off completely if required.

The equations~\cref{eq:advection} have a unique smooth solution $c\in H^1(\Omega_T)$ provided the data in~\cref{eq:advection}, namely $\bm{u}$, $e$, $c_0$, $e_0$, $g$, $g_0$, $g_\partial$ and $c_\partial$, $e_\partial$ are smooth enough, i.e. of $H^1(\Omega_T)$ class, and the initial condition agrees with boundary conditions (see~\cite[p.484]{Quarteroni2008_NumApproxPDEs}). For the case of less regular solutions, namely just measurable data, the equations~\cref{eq:advection} still possess the unique weak solution of $L^2(\Omega_T)$ class (see~\cite[p.220]{Bardos1970}), and the discontinuities/steep gradients are propagated forward by the characteristics. For the sake of completeness, we note that the solutions of~\cref{eq:advection} can receive an even more general representation, namely that of solutions in the space of measures~\cite{Raviart1985} (for the case $\Omega=\mathbb{R}^n$). Since our goal is not in proving the convergence of the proposed method for the most general case, but rather in demonstrating how one can estimate the state of~\cref{eq:advection} efficiently without compromising much the estimation precision, in what follows, we assume that the data are smooth enough so that $c$ is at least continuous in $\bm{x}$ and has the classical gradient of $L^2$-class.

\paragraph{Observation equation} 
We further assume that a function $y$ is observed:
\begin{equation}
  \label{eq:obs}
  y(\bm{x}_j,t) = c(\bm{x}_j,t)+\eta_j(t)\,, \quad j=1\dots N_s\,,
\end{equation}
where $\bm{x}_j\in\Omega$ denotes the position of a sensor, and a network of $N_s$-sensors is deployed in $\Omega$.
Finally, the observation noise $\bm{\eta}:=(\eta_1\dots\eta_{N_s})^\top$ is modelled as a realisation of a vector-valued random process with zero mean and uncertain but bounded correlation function:
\begin{equation}
  \label{eq:ellips_Y}
  \int_0^T E(\bm{\eta}(t) \cdot R^{-1}(t) \bm{\eta}(t))dt\le 1\,,
\end{equation}
where $R(t)$ is a given symmetric positive definite continuous weighting matrix with continuous inverse. This assumption reflects the fact that the second moments of $\bm{\eta}$ are known up to a bounded error, e.g. when the moments are estimated by an empirical moment estimator, as often happens in practice.

\emph{Problem statement:} given state equation~\cref{eq:advection} together with observations~\cref{eq:obs} and the uncertainty descriptions~\cref{eq:model_error},~\cref{eq:ellips_Y} our goal is to design an efficient numerical algorithm estimating the quantity $c(\bm{x},t)$ from $(y(\bm{x}_1,t)\dots y(\bm{x}_{N_s},t))^\top$.

\section{Mathematical preliminaries}
In this section we provide well-known information on DG discretization for advection equations, and minimax filtering for linear ODEs. This material forms a basis for~\cref{sec:main}. 
\subsection{Discontinuous Galerkin formulation}
\label{sec:DGform}
Applying the DG method to equation~\cref{eq:advection} is along the same lines as its application to a general non-linear conservation law. For more details on what follows, see~\cite{hesthaven2008nodal}. First, the domain, $\Omega$ is divided into $K$ of non-overlapping elements, $D_k$, i.e., $\Omega 	\simeq \Omega_h = \bigcup\limits_{k=1}^{K} D^k$, where we choose the elements, $D^k$ to be rectangular. The restriction of the state $c$ onto the element $D_k$ is denoted by $c^k$. The latter can be approximated by $c_h^k$, which is expressed as the series,
\begin{equation}
\label{eq:seriesDG}
c_h^k(\bm{x},t)= \sum_{i=1}^{N+1} c^k(\bm{x}_i^k,t) \ell_i^k(\bm{x}), \quad \bm{x} \in D^k,
\end{equation}
where $\ell_i^k(\bm{x})$ are Lagrange interpolating polynomials in two dimensions defined on points $\bm{x}_i^k$. These points are taken to be quadrature points for Legendre polynomials, specifically Legendre-Gauss-Lobatto (LGL) points. The series~\cref{eq:seriesDG} is a nodal expansion representing $c_h^k$, i.e. $c_h^k(\bm{x}_i^k,t) = c^k(\bm{x}_i^k,t)$. For more details on this, see~\cite{hesthaven2008nodal}. We define $e_h^k$ and $e_k^0$, $e_\partial^h$, $g_h^k$, $g_{h,0}^k$ and $g_{h,\partial}^k$ analogously. Substituting $c_h^k$ into~\cref{eq:advection}, we form the residual $R_h^k$ on a single element:
\begin{equation}
\label{eq:residual}
R_h^k = \partial_t c_h^k(\bm{x},t)+ \bm{u} \cdot \nabla c_h^k(\bm{x},t) - g_h^k(\bm{x},t) e_h^k(\bm{x},t), \quad \bm{x} \in D^k.
\end{equation}
The Galerkin method involves taking $\ell_i^k$ as test functions (i.e. same as the expansion/trial functions) and forcing the residual to be orthogonal to each of these test functions. Doing this, and then using integration by parts \[
\int_{D^k}  \left( (\bm{u}\cdot \nabla c_h^k) \ell_n^k + \nabla \cdot \bm{u} c_h^k\ell_n^k + (c_h^k\bm{u})\cdot \nabla \ell_n^k \right)  d\bm{x} = \int_{\partial D^k} (\bm{\hat{n}}\cdot \bm{u}) \ell_n^k c_h^k d\bm{x}\,,
\]
to move the spatial derivatives off the state $c$ and onto the test functions gives\footnote{Recall that $\nabla \cdot \bm{u}=0$.} the following weak statement on the element $D^k$ ($n=1 \ldots N+1$):
\begin{equation}
\label{eq:weakStatement}
\begin{split}
\int_{D^k}  \left( \partial_t c_h^k \ell_n^k - \bm{f}_h^k\cdot \nabla \ell_n^k\right)  d\bm{x} = -\int_{\partial D^k} \bm{\hat{n}}\cdot \bm{f}^* \ell_n^k d\bm{x} + \int_{D^k} g_h^k e^k_h \ell_n^k d\bm{x} \,,
\end{split}
\end{equation}
where $\bm{\hat{n}}(\bm{x})=(\hat n_x,\hat n_y)^\top$ is the unit outward vector which is normal to $\partial D_k$ at $\bm{x}$, $\bm{f}_h^k=(uc_h^k,vc_h^k)^\top$ and $\bm{f}^*$ is the numerical flux function which we take to be the local Lax-Friedrichs flux:
\begin{equation}
\label{eq:LF}
\bm{f}^*(c_i,c_e,\bm{u}_i,\bm{u}_e) = \frac{c_i \bm{u}_i + c_e \bm{u}_e }{2}+\frac{\mu}{2}\bm{\hat{n}}(c_i-c_e),
\end{equation}
where subscripts $i$ and $e$ refer respectively to the interior and exterior values at a point on the boundary, and $\mu$ is the maximum absolute value of the signal speed normal to the boundary at that point, i.e.,
\begin{equation}
\label{eq:LFconst}
\mu=\max\{| \bm{u}_i(\bm{x})\cdot\hat{\bm{n}}|, | \bm{u}_e(\bm{x})\cdot\hat{\bm{n}}| \}
\end{equation}
The surface integral in~\cref{eq:weakStatement} allows the elements to `communicate' with one another by imposing the values for $c_h^k$ at the boundary of $D^k$ from the adjacent elements. In fact, the boundary values are imposed in a weak sense~\cite[p.483]{Quarteroni2008_NumApproxPDEs} as opposed to the strong/classical sense when the solution takes exactly the prescribed value at the boundary~\cite[p.482]{Quarteroni2008_NumApproxPDEs}. This strategy of weak boundary conditions agrees well with our problem statement: indeed, we assume that boundary conditions in~\cref{eq:advection} are given up to an uncertain function, and so it does make sense to impose boundary conditions ``on average'' or in the weak sense.
Since we are using rectangular elements, the surface integral is just the sum of four line integrals, each one over one face of the element $D^k$. The exterior solution values $c_e$ in \cref{eq:LF} and \cref{eq:LFconst} refer to the value of $c^k_h$ at $\bm{x}$ at the boundary of a neighbouring element that shares the boundary with $D^k$. If $\bm{x}$ belongs to the physical boundary of the domain, $\Gamma$, the exterior value is determined by the third equation of~\cref{eq:advection}. More specifically, if the flow direction at $\bm{x}$ is `into' the domain, i.e. $\bm{\hat{n}}(\bm{x})\cdot \bm{u}(\bm{x})<0$, then $c_e$ is set to either the value of $c^k_h$ at $\bm{x}$ at the boundary of a neighbouring element or the value prescribed by the 2nd equation of~\cref{eq:advection}. If, on the other hand, the flow direction at that boundary point is `out of' the domain (i.e. $\bm{\hat{n}}(\bm{x})\cdot \bm{u}(\bm{x})>0$), then a free exit boundary condition is imposed at that point by setting $c_e$ equal to the interior value $c_i$. Note that the Lax-Friedrichs numerical flux adds artificial viscosity at the element interfaces in its handling of the jump in state, $c_i-c_e$. This smooths discontinuitues that may occur at those interfaces.

The weak DG formulation of~\cref{eq:advection} on a single element $D^k$ can be written as:
\begin{equation}
\label{eq:weakDG}
\begin{split}
M^k \dfrac{d \bm{c}_h^k}{dt} - S_x^\top \bm{f}_x - S_y^\top \bm{f}_y =-\sum\limits_{i=1}^4 (-1)^i M_e^{k,i} \bm{f}^*_i+M^kG^k(t)\bm{e}^k(t)\,,\quad \\
\bm{c}_h^k(0) = \bm{c}_0^k + G_{0}^k\bm{e}_0^k\,,
\end{split}
\end{equation}
where $M^k$, $S_x$ and $S_y$ are the mass and stiffness matrices with the latter corresponding to advection in the $x$- and $y$-directions. The vectors $\bm{c}_h^k$, $\bm{c}_0^k$, $\bm{e}^k$ and $\bm{e}_0^k$ are grid-functions representing $c_h^k$, $c_0$, $e_0$, $e$ respectively on the $(N+1)^2$ LGL quadrature points of the element $D^k$, and the vectors $\bm{f}_x$ and $\bm{f}_y$ are grid functions representing the first and second components of $\bm{f}_h^k$ respectively. $G^k$ and $G^k_0$ are diagonal matrices with the grid-functions $\bm{g}_h^k$ and $\bm{g}_{h,0}^k$ on the diagonal respectively. The matrices $M_e^{k,i}$ are \emph{edge}-mass matrices for the element $D^k$ on face $i$, where the faces, $i= 1 \ldots 4$, are ordered: \emph{left}, \emph{right}, \emph{lower}, \emph{upper}. These matrices act on the vector $\bm{f}^*_i$ which is a grid function representing the numerical flux,~\cref{eq:LF}, over each node on face $i$.

As noted above, the weak DG formulation of~\cref{eq:advection} on a single element $D^k$ leads to the ODE~\cref{eq:weakDG}, which describes the time evolution of the vector $\bm{c}_h^k$, the grid-function representing $c_h^k$ on the $(N+1)^2$ LGL quadrature points on element $D^k$. The source term $\bm{f}^*_i$ of~\cref{eq:weakDG} depends upon vectors $\bm{c}_h^s$, $s\ne k$ which approximate $c_h^s$ on $D^s$. It is responsible for the ``exchange of information'' between the ``local'' approximations of $c_h^s$. This ``communication mechanism'' plays an important role in the distributed filtering algorithm presented in~\cref{sec:localMF}.

\subsection{Minimax Filter}
\label{sec:MF}
Let $\bm{x}(t)\in\mathbb{R}^n$ be a state vector of the following equation: 
\begin{equation}
\label{eq:state}
\dfrac{d \bm{x}}{dt} = A(t) \bm{x} + \bm{b}(t) + V\bm{e}(t)\,, \bm{x}(0) = \bm{x}_0+G_0\bm{e}_0,
\end{equation}
where $A,V, G_0$ are given matrices, $\bm{b}\in\mathbb{R}^n$ is a given vector-function representing a source term, $\bm{e}\in\mathbb{R}^m$ is a measurable function representing the model error, and $\bm{e}_0\in\mathbb{R}^p$ is an unknown vector representing the error in the initial condition $\bm{x}_0\in\mathbb{R}^n$. We assume that $\bm{e}_0$ and $\bm{e}$ are uncertain but bounded, and belong to the following bounding set: 
\begin{equation}
  \label{eq:ellips}
E_T:=\{(\bm{e}_0,\bm{e}): \bm{e}_0\cdot Q_0 \bm{e}_0 + \int_{0}^{T} (\bm{e}(t)\cdot Q(t) \bm{e}(t) dt \le 1\}
\end{equation}
provided $Q_0$ and $Q(t)$ are given symmetric positive definite matrices, and $Q(t)$ has bounded (in $t$) inverse. It is not hard to see that $E_T$ is an ellipsoid in the space $\mathbb{R}^p\times L^2(0,T,\mathbb{R}^m)$.\\
We further assume that a vector-function $\bm{y} \in\mathbb{R}^s$ is observed: 
\[
\bm{y}(t) = H(t) \bm{x}(t) + \bm{\eta}(t)\,,
\]
where the assumptions on the observation noise $\bm{\eta}\in\mathbb{R}^s$ are the same as in~\cref{sec:probStat}. 

Let us introduce the minimax filter equations. The minimax filter, $\hat{\bm{x}}$ solves the following system: 
\begin{align}
\frac{dP}{dt} &= A P + P A^\top + VQ^{-1}V^\top - P H^\top R^{-1} H P \,, P(0) = G_0Q_0^{-1} G_0^\top ,\\
\frac{d \hat{\bm{x}}}{dt} &= A{\hat{\bm{x}}} + \bm{b}(t) + P H^\top R^{-1}({\bm{y}} - H{\hat{\bm{x}}})\,, \hat{\bm{x}}(0) = \bm{x}_0
\end{align}
It can be shown that the worst-case state estimation error associated with the minimax filter $\hat{\bm{x}}$ can be expressed as follows: 
\begin{equation}
\max_{(\bm{e}_0,\bm{e})\in E_T, R } E(\bm{x}_j(t)- \hat{\bm{x}}_j(t))^2 = P_{jj}(t).
\end{equation}
It is important to note that the worst-case mean-squared error of the minimax estimate is optimal in that any other estimate of the state that is either linear or non-linear in the observations would have a corresponding worst-case mean-squared error larger than that of the minimax~\cite{SZAPOGN_TAC16}.   

\section{Main result}
\label{sec:main}

As noted in the introduction, our approach relies upon the ``discretize-then-optimize'' strategy composed of the following steps:
\begin{itemize}
\item [1)] the domain, $\Omega$ is divided into $K$ of non-overlapping elements, $D_k$, and the weak formulation of the DG method is applied to~\cref{eq:advection} on $D^k$, namely~\cref{eq:advection} is substituted by a system of ODEs which describe the evolution (in time) of a vector approximating $C$ on an element $D^k$ (\cref{sec:DGform});
\item [2)] for each element $D^k$, the minimax state estimator is applied to estimate the state of the system of ODEs by using the data localised in $D^k$; the estimators, located at adjacent elements, `communicate' with one another by means of the boundary integral interconnection mechanism of the DG method (\cref{sec:localMF}).
\end{itemize}

\subsection{Distributed minimax filter}
\label{sec:localMF}

To proceed we first transform~\cref{eq:weakDG} derived in~\cref{sec:DGform} into a form suitable for the application of the minimax state estimator described in~\cref{sec:MF}. We introduce an affine transformation:
\begin{equation}
\begin{split}
\label{eq:affine}
A^k(t) \bm{c}_h^k(t) + \bm{b}^k(t,\bm{c}_h^s) + W^k\bm{e}^k_\partial(t) := (M^k)^{-1} (S_x^\top \bm{f}_x(t) + S_y^\top \bm{f}_y(t)) \\
-(M^k)^{-1}\sum\limits_{i=1}^4 (-1)^iM_e^{k,i} \bm{f}^*_i(t)
\end{split}
\end{equation}
where the first term on the right hand side accounts for the advection of the state quantity, $\bm{c}_h^k$ within the element $D^k$ in the absence of boundary effects, and the second term accounts for inter-element boundary fluxes. Of those two terms, the first is linear in $\bm{c}_h^k$, as the mass and stiffness matrices are independent of $\bm{c}_h^k$, while we see from the definition of $\bm{f}_x$ and $\bm{f}_y$ in~\cref{sec:DGform} that they depend linearly on $\bm{c}_h^k$. The second term, however is affine in $\bm{c}_h^k$, which we can see if we look at the numerical flux given by~\cref{eq:LF}, which is approximated by the grid-function, $\bm{f}^*_i$ above. The interior state, $c_i$, in that equation is approximated by $\bm{c}_h^k$ restricted to the nodes on the appropriate boundary. However, the other terms in that equation do not depend on $\bm{c}_h^k$ so are treated as an external source term. As a result, in addition to containing the stiffness terms, $A^k \bm{c}_h^k$ also absorbs part of the flux term, while $\bm{b}^k$ represents only the component of the flux that depends on the state $\bm{c}_h^s $ on elements $D^s$ neighbouring $D^k$. Finally, $\bm{e}^k_\partial $ represents the component of the flux that depends on $e_\partial$ which is only the case when $D^k\cap\Gamma^{in}\ne\varnothing$. The term $\bm{e}^k_{\partial}$ is in effect a flux term induced by the error at the physical boundary, i.e., $\bm{e}^k_{\partial}=-(M^k)^{-1}\sum\limits_{i=1}^4 (-1)^iM_e^{k,i} \bm{f}^*_{\partial}(t)$ where $\bm{f}^*_{\partial}(t)$ is a grid-function approximating the Lax-Friedrichs flux,~\cref{eq:LF}, with external state, $c_e = e_{\partial}$. The matrix, $W^k$ restricts $\bm{e}_{\partial}$ to the boundary $D^k$ when $D^k$ has a boundary that intersects with $\Gamma$, otherwise $W^k$ is zero. Note that the effect of $c_\partial$ is absorbed by $\bm{b}^k$. In other words, $\bm{b}^k$ represents the ``known'' part of the boundary conditions (either from $\partial D^k$ or from $\Gamma^{in}$), and $\bm{e}^k_\partial $ accounts for the uncertain part of the boundary conditions coming from $e_\partial$. In what follows, we will refer to $A^k$ as the state transition matrix, or system matrix for short. Having said this, we can rewrite~\cref{eq:weakDG} as follows:
\begin{equation}
\label{eq:affineSys}
\dfrac{d \bm{c}_h^k}{dt} = A^k(t) \bm{c}_h^k + \bm{b}^k(t,\bm{c}_h^s) + V^k\bm{\tilde e}^k(t)\,, \bm{c}_h^k(0) = \bm{c}^k_0+G_0^k\bm{e}^k_0\,.
\end{equation}
where $V^k:=(G^k\,\,\, W^k)$ and $\bm{\tilde e}^k=(\bm{e}^k,\bm{e}^k_\partial)^\top$. Note that $G^k=0$, $W^k=0$ provided $g=0$ and $g_\partial=0$ respectively.

Now, let us introduce a bounding set for the uncertain vectors $\bm{e}^k_0$ and $\bm{\tilde e}^k$. We recall that $\bm{e}^k_0=(c_0(\bm{x}_1^k) \dots c_0(\bm{x}_{(N+1)^2}^k))^\top$, and $\bm{e}^k$, $\bm{e}^k_\partial$ are defined analogously. Define $Q_0^k:=\operatorname{diag}(q_0(\bm{x}^k_1)\dots q_0(\bm{x}^k_{(N+1)^2}))$, the restriction of $q_0$ onto the LGL grid of $D^k$. Let $Q^k(t)$ denote the restriction of $q$ onto the LGL grid, defined in the same way. Let $Q^k_\partial(t)$ denote the restriction of $q_\partial$ onto the LGL grid points located at the boundary of $D^k$. Since $q_0,q_\partial, q>q^\star>0$, it follows that $Q_0^k$, $Q_\partial^k$ and $Q^k$ are positive definite. By~\cref{eq:model_error}, we get that
\begin{equation}
  \label{eq:ellips_Dk}
\bm{e}^k_0\cdot Q^k_0 \bm{e}^k_0 + \int_{t}^{t+s} (\bm{\tilde e}^k(\tau)\cdot \tilde Q^k(\tau) \bm{\tilde e}^k(\tau) dt \le (1+2s) (N+1)^2\,,
\end{equation}
provided $\tilde Q^k:=\operatorname{diag}(Q^k,Q_\partial)$. This approximation represents an ellipsoid containing $\mathcal{C}^k$, the restriction of the hypercube defined by~\cref{eq:model_error} onto the LGL grid within $D^k$. Clearly, the ellipsoid ``contains more uncertainty'' than the restricted hypercube $\mathcal{C}^k$: indeed, this is indicated by the presence of the factor $(N+1)^2$. However, in practice, $N$ (the degree of Lagrange polynomials) is typically not taken to be higher than the low integers. The numerical precision of the DG method is increased by fixing $N$ and refining the partition of $\Omega$ by introducing more elements $D^k$ of smaller area. For example, we use $N=3$ and $70\times70$ rectangular partition of $\Omega$ to advect a satellite image of the cloud optical depth (see~\cref{sec:COD}). Therefore, the factor $(N+1)^2$ has a positive effect (at least for $N\le 5$) as it allows one to include uncertain source terms $\bm{\tilde e}^k$ with slightly larger energy; for instance, setting $\bm{\tilde e}^k = \bm{\tilde e}^k_1 + \bm{d}^k$ where the latter term accounts for discretization error (provided it is small enough, which in turn can be achieved by taking a large enough number of small elements)  making the resulting state estimate more robust by increasing the worst-case estimation error (see~\cref{p:1}). We refer the reader to~\cite{ZhukSISC13} for further discussion on including the discretization error into the ellipsoid. 

Finally, we recall that a network of $N_s$-sensors is deployed in $\Omega$. We denote by $D^k_{obs}:=\{\bm{x}_{j_1} \ldots \bm{x}_{j_{M_k}}\}$ the locations of the sensors that belong to $D^k$ and define \[
\bm{y}^k=(y(\bm{x}_{j_1},t)\ldots y(\bm{x}_{j_{M_k}},t))^\top\,\quad  H^k=\{\ell^k_n(\bm{x}_{j_{M_k}})\}_{n,j=1}^{M_k,(N+1)^2},
\] the restriction of $y$ defined by~\cref{eq:obs} onto $D^k_{obs}$, and the interpolation matrix $H^k$ mapping LGL grid to $D^k_{obs}$. Similarly, we define $\bm{\eta}^k=(\eta_{j_1}\dots \eta_{j_{M_k}})^\top$, the restriction of the observation noise onto the element $D^k$. Let $\pi$ denote a matrix of norm 1 mapping $\bm{\eta}$ to $\bm{\eta}^k$, and define $R^k:=\pi^\top R \pi$. If $D^k$ does not contain a single sensor we set\footnote{It will become apparent after Proposition~\cref{p:1} that setting $H^k=0$ is mathematically equivalent to the ``no observation'' case: indeed, the state estimator is coupled to $\bm{y}^k$ by means of $H^k$ so this coupling becomes trivial if $H^k=0$ and the state estimator reduced to the state equation~\cref{eq:affineSys} with no uncertainty, e.g. $\bm{e}_\partial=0$, $\bm{e}^k_0=0$ and $\bm{e}^k=0$.} $H^k:=0$. As a result, the observations take the following form: \[
\bm{y}^k = H^k \bm{c}_h^k(t) + \bm{\eta}^k\,,
\]
Now, following~\cite{SZAPOGN_TAC16} we state the following proposition:
\begin{proposition}\label{p:1}
  Assume that $\bm{c}^k$ solves~\cref{eq:affineSys} and the uncertain parameters are bounded according to~\cref{eq:ellips_Dk}. Given observations $\bm{y}^k$ we define the minimax estimate $\hat{\bm{c}}^k_h$ on $D^k$ as follows:
\begin{align}
\frac{dP^k}{dt} &= A^k P^k + P^k (A^k)^\top + V^k(\tilde{Q}^k)^{-1}(V^k)^\top\label{eq:ric1}\\
& - P^k (H^k)^\top (R^k)^{-1} H^k P^k \,, \qquad P^k(t) = P^k_{prev}(t)\,,\label{eq:ric2}\\
\frac{d \hat{\bm{c}}^k_h}{dt} &= A^k{\hat{\bm{c}}_h^{k}} + \bm{b}^k(t,\hat{\bm{c}}^{s}_h)\label{eq:lMF1}\\
&  + P^k (H^k)^\top (R^k)^{-1}({\bm{y}}^k - H^k{\hat{\bm{ c}}^k_h})\,, \qquad \hat{\bm{c}}^k_h(0) = \bm{c}^k_0 \label{eq:lMF2}
\end{align}
where $\hat{\bm{c}}^{s}_h$ is the minimax estimate on elements $D^s$ neighbouring $D^k$. It then follows that
\begin{equation}
  \label{eq:lMF_err}
\max_{\tilde{\bm{e}}^k,\bm{e}^k_0,E \bm{\eta}^k (\bm{\eta}^k)^\top } E(c_h^k(\bm{x}_j,t+s)- \hat c_h^k(\bm{x}_j,t+s))^2 = (1+2s)(N+1)^2P_{jj}(t+s)\,,
\end{equation}
i.e. the worst-case mean-squared error of the minimax estimate $\hat c_h^k(\bm{x}_j,t+\delta t)$ of the value $c_h^k(\bm{x}_j,t+\delta t)$ is given by the $j$th diagonal element of the unique symmetric positive definite solution of the Riccati equation~\cref{eq:ric1}-\cref{eq:ric2}.
\end{proposition}
The proof of the proposition is given in~\cref{s:proofs}. Note that \cref{p:1} reflects the sequential nature of the estimation process: at $t=0$ we initialize the process by setting $P^k_{prev}(0):= G_0^k(Q_0^{k})^{-1}G_0^k$, which describes the a priori bound for the initial condition error, and compute $P^k$ and $\hat{\bm{c}}^k_h$ on $(0,s)$ by solving~\cref{eq:ric1}-\cref{eq:lMF2}; $t$ is then set to $s$, $P_{prev}(t)$ is set to $P^k(t)$ and~\cref{eq:ric1}-\cref{eq:lMF2} is solved again on $(t,t+s)$. The factor $(1+2s)(N+1)^2$ in~\cref{eq:lMF_err} comes from the ellipsoidal approximation of the discretized hypercubes~\cref{eq:model_error}, which is required to formulate the minimax filter~\cite{Zhuk2009d}. As a result, the mean-squared worst-case estimation error is inflated after each estimation step at $t+s$ by the constant factor $(1+2s)(N+1)^2$ so that at time $T$ the estimation error is given by $(1+2T)(N+1)^2 P_{jj}$. 
We stress that the equation for $P^k$ is composed of the following parts:
\begin{itemize}
\item the linear part~\cref{eq:ric1} represents a Lyapunov operator that describes dynamics of the estimation error in the absence of observations (e.g. $H^k=0$),
\item the nonlinear part~\cref{eq:ric2} represents the reduction in the estimation error due to observations.
\end{itemize}
Similarly, the state estimator equation consists of the model~\cref{eq:lMF1} and the innovation part~\cref{eq:lMF2}. Both ~\cref{eq:ric2} and \cref{eq:lMF2} disappear if $D^k$ does not contain a single sensor.

We stress that the gain $P^k$ does not depend explicitly on the gains $P^s$ at the neighbouring elements. In contrast, the minimax estimates $\hat{\bm{c}}^k_h$ at elements $D^k$ are advanced independently over the time window $[t,t+s]$, and the communication terms $\bm{b}^k(t,\hat{\bm{c}}^{s}_h)$ are then updated at $t+s$. In this way, the filters on elements with no observations receive information from the elements with observations, so that the localised observations are, in fact, spread around the entire domain by the communication terms $\bm{b}^k(t,\hat{\bm{c}}^{s}_h)$. This same mechanism provides an implicit communication between $P_k$: indeed, $\bm{b^k}$ depends on $\hat{\bm{c}}^{s}_h$ at time $t_s$, and this changes $\hat{\bm{c}}^k_h$ at $t_{s+1}$, which is in turn reflected in the local system matrix $A^k$. The latter, in turn, modifies $P^k$.

\subsection{Time discretisation}
\label{sec:timeDisc}
It was noted in~\cite{ZhukSISC13} that equations~\cref{eq:ric1}-\cref{eq:lMF2} must be discretized by a method (e.g. symplectic Runge-Kutta method of order $p$) preserving quadratic invariants of the estimation error dynamics, i.e. the discrete estimate should verify the equality~\cref{eq:lMF_err}. The latter holds true for the symplectic Mobius integrator proposed in~\cite{ZhukCDC14} to solve the matrix differential Riccati equation~\cref{eq:ric1}. The basic idea behind the Mobius transformation is to make use of the fact that the solution of the Riccati equation induces a flow on a Grassmannian manifold. This flow is called a Mobius transformation. It may be constructed by solving an associated linear Hamiltonian system: indeed, the solution of the Riccati equation, $P^k$, can be expressed in the form $P^k=V^k(U^k)^{-1}$ provided
\begin{equation}
\label{eq:UV}
\begin{pmatrix}
  \dot U^k\\ \dot V^k  
\end{pmatrix}
 = 
\begin{pmatrix}
-A^{k \top} &  (H^k)^\top(R^k)^{-1}H^k \\ \bar Q_k & A^k
\end{pmatrix}\begin{pmatrix}
U^k(t) \\ V^k(t)
\end{pmatrix}
, \quad
\begin{pmatrix}
U^k(t) \\ V^k(t)
\end{pmatrix}
=
\begin{pmatrix}
I \\ P^k(t)
\end{pmatrix}
\end{equation}
where $\bar Q_k:=(1+2s)(N+1)^2V^k\tilde (Q^k)^{-1}(V^k)^\top$, and the initial gain, $P^k(0):= (1+2s)(N+1)^2G_0^k(Q_0^{k})^{-1}G_0^k$. The Hamiltonian system~\cref{eq:UV} can be solved by using symplectic Runge-Kutta methods of order $2$ and thus avoid numerical instabilities associated with Mobius transform by means of a reinitialization: namely the gain at time level $j+1$ is given by $P^k_{j+1}=V_{j+1}^k (U_{j+1}^k)^{-1}$ provided
\begin{equation}\label{eq:UVj}
  \begin{split}
  \begin{pmatrix}
    I+\frac{\Delta T}2 (A_j^k)^\top & -\frac{\Delta T}2(H^k)^\top(R^k_{j+1})^{-1}H^k  \\ -\frac{\Delta T}2 (\bar{Q}^k)^{-1} & I-\frac{\Delta T}2 A_j^k
  \end{pmatrix}
\begin{pmatrix}
  U_{j+1}^k\\V_{j+1}^k
\end{pmatrix}
\\=
  \begin{pmatrix}
    I-\frac{\Delta T}2 (A_{j}^k)^\top & \frac{\Delta T}2  (H^k)^\top(R^k_j)^{-1}H^k \\ \frac{\Delta T}2  (\bar{Q}^k)^{-1}& I+\frac{\Delta T}2 A_{j}^k
  \end{pmatrix}
  \begin{pmatrix}
  I\\P_j^k
  \end{pmatrix}
  \end{split}
\end{equation}
where $\Delta T$ is the time-step, and $s:=\Delta T$. The reinitialization is in that $U^k_j$ is set to the identity matrix and $V^k_j=P_j^k$ so that $(U^k_{j+1})^{-1}$ is well-conditioned. The resulting symplectic Mobius integrator~\cref{eq:UVj} is stable and preserves symmetry and positivity of the Riccati matrix. Also it preserves all quadratic invariants of the estimation error dynamics including~\cref{eq:lMF_err}. Finally, the filter equation in~\cref{eq:lMF1}-\cref{eq:lMF2} is also solved using the implicit midpoint method, giving:
\begin{equation}
  \label{eq:filter}
  \begin{split}
    &\left(I-\frac{\Delta T}2 A_{j+1}^k+\frac{\Delta T}2 P_{j+1}^k (H^k)^\top(R^k_{j+1})^{-1}H^k \right) (\hat{\bm{c}}_h^k)_{j+1}= \Delta T\bm{b}^k_j\\
    &\left(I+\frac{\Delta T}2 A_{j}^k-\frac{\Delta T}2 P_{j}^k (H^k)^\top(R^k_{j})^{-1}H^k \right)(\hat{\bm{c}}_h^k)_j \\
    &+\frac{\Delta T}{2}\left(\frac{P_{j+1}^k (H^k)^\top(R^k_{j+1})^{-1}\bm{y}_{j+1}^k +P_j^k  (H^k)^\top(R^k_{j})^{-1} \bm{y}_j^k}{2}\right)\,,
  \end{split}
\end{equation}
Note that the above scheme is explicit in $\bm{b}^k_j$ as at time level $j$ we only have $\bm{b}^k_j$.

\subsection{Varying trust in observations}
\label{sec:varyingTrust}

When filtering, the amount of trust placed in the observations, $\bm{y}^k$, is regulated by the symmetric positive definite matrix, $R^k$: small eigenvalues of $R$ represent high trust and the reverse. Intuitively, low trust (high eigenvalues) reduces the ``rate'' at which local filter assimilates the observations by reducing the impact of the innovation term~\eqref{eq:lMF2}. This simple fact is used to assimilate sparse in time observations in a stable fashion. Indeed, if the observations are available at discrete time instants, i.e. the matrix $H^k$ switches between $0$ and identity matrix, the numerical scheme~\eqref{eq:lMF1}-~\eqref{eq:lMF2} could quickly become unstable due to the hyperbolic nature of the problem~\eqref{eq:advection}. To overcome this, we suggest the following procedure: instead of switching $H^k$ between $0$ and $1$ to mimic the presence/absence of the observations, we fix $H^k$ and vary the trust in $\bm{y}^k$ by multiplying $R^k$ by a scalar $r^k>0$ which dictates the trust placed in the observations associated with element $k$. Small $r^k$ indicates high trust, while large $r^k$ indicates low trust. More specifically, the algorithm assumes that no data are available at $t=0$, so $r^k$ is initialised to a high value. When observations become available, say at $t=t_1$, then $r^k$ should be decreased from its default high value to a low value in order to increase the trust in the observations for element $k$. However, sharply decreasing $r^k$ can cause the system to become numerically unstable, so instead of doing this, we reduce its value over time, thus increasing the trust gradually and avoiding instabilities. The observations should however be assimilated relatively fast to avoid lag, so reducing $r^k$ at each computational time-step would not be ideal unless $\Delta T$ was set to a small value (i.e. smaller than necessary to ensure numerical stability of the scheme). This, however, would be computationally inefficient so instead, we introduce a second, smaller time-step, $\Delta T_s$, where $n_s  \Delta T_s = \Delta T$, so that when data become available, we can switch to using the smaller time-step, $\Delta T_s$, running the filter $n_s$ times over $[t_1, t_1 +  \Delta T]$ while varying $r^k$ over that time-interval. The trust at the beginning of the time interval is small (high $r^k$) and should also be small at the end. Thus, we vary $r^k$ in such a way that it attains a specific minimum value associated with high trust mid-way through the time interval. We achieve this by reducing $r^k$ for the first $n_s/2$ small time-steps (setting $n_s$ to be even) and then increasing it back to a high value over the remaining $n_s/2$ steps. This procedure resembles the (weak) approximation of the Dirac delta-function by Gaussian densities.

The distributed filtering algorithm is summarised as~\cref{alg:localFilter}, where we use the following notation:
\\
\begin{itemize}
\item $\Delta T$: Standard computational time-step
\item $N_t$: Number of standard time-steps
\item $N_c$: Number of time-steps used at current time-level
\item $n_s$: Number of small time-steps for filtering (even)
\item $K$: Number of DG elements
\item $\Delta T_c$: Time-step at current time-level
\item $E_s$: Set of DG elements on which observations may be available
\item $r^k$: Observation trust parameter
\item $r_{tr}$: High trust value of $r^k$
\item $r_d$: Low trust value of $r^k$
\item $\tau$: Constant for varying trust, where $r_d\tau^{-{n_s}/2} = r_{tr}$
\end{itemize}
\begin{algorithm}
\caption{Distributed filtering for 2D advection}
\label{alg:localFilter}
\begin{algorithmic}
\STATE $r^k \gets r_d$
\FOR{$i=1$ to $N_t$}
    \STATE $t \gets (i-1)\Delta T$
    \IF{Observations are available on any element at time, $t$}
	\STATE $N_c \gets n_s$
	\STATE $\Delta T_c \gets \Delta T / n_s$
    \ELSE
	\STATE $N_c \gets 1$
	\STATE $\Delta T_c \gets \Delta T$
    \ENDIF

    \FOR{$j=1$ to $N_c$}
	\STATE $t \gets t + j\Delta T_c$
	\STATE Update boundary conditions and advection field for time $t$
        \FOR{$k=1$ to $K$}
            \STATE Compute elemental system matrix, $A^k$, and vector, $\bm{b}^k$
	    \IF{Element $k \in E_s$ }
                \STATE \textit{\textbf{Filter:}} Compute gain, $P^k$, and estimate, $\hat{\bm{c}}_h^k(t)$, using~\cref{eq:UVj} and~\cref{eq:filter} with current time-step, $\Delta T_c$
		\IF{$N_c = n_s$}
		    \STATE \textit{\textbf{Vary trust in observations:}}
		    \IF{$j \leq n_s/2$}
		        \STATE $r^k \gets r^k /\tau$
		    \ELSE
			\STATE $r^k \gets r^k \tau$
		    \ENDIF
		\ENDIF
            \ELSE
                \STATE \textit{\textbf{Solve forward:}} Set $H^k=0$ and compute $P^k$ and $\hat{\bm{c}}_h^k$, using~\cref{eq:UVj} and~\cref{eq:filter} 
            \ENDIF
        \ENDFOR
        \STATE Assemble full estimate $\hat{\bm{c}}_h(t)$ at time $t$ by combining all $K$ elemental estimates $\hat{\bm{c}}_h^k(t)$
    \ENDFOR
    
\ENDFOR
\end{algorithmic}
\end{algorithm}
Note that there are two nested time-loops in~\cref{alg:localFilter}. For the outer loop (using control variable, $i$) time is advanced using the standard computational time-step, $\Delta T$, and for the inner one (using control variable, $j$) time is either incremented by the standard computational time-step for a single iteration in the case that no observations are available, or by the small time-step, $\Delta T_s$ for $n_s$ iterations when observations are available on any element. The time, $t$, that is set inside the outer loop is the time at the beginning of the current time-level. It is at this time that a check is performed for available observations. Inside the inner time-loop, $t$ is then updated to be the estimate/solution time. This will be either $\Delta T$ or $\Delta T_s$ beyond the time at the beginning of the current time-level, depending on whether or not observations are available.

Note also that the loop over DG elements (using control variable, $k$) can be easily parallelised as it is called within the inner time-loop and hence runs over a single time-level at a time. As a result, the elemental estimates/solutions can be obtained in any order, since~\cref{eq:UVj} and~\cref{eq:filter} form a two-level scheme whereby the solution at the current time-level is obtained using the solution at the previous level.

\section{Numerical Experiments}
\label{sec:experiments}

In this section, we describe numerical experiments in both synthetic and real scenarios. The real scenario is that of cloud motion where we use a sequence of velocity fields that are obtained using an optical flow estimation procedure.  The observations used are satellite images depicting cloud optical depth. We employ the distributed filter here on a high-resolution grid on a domain with constant inflow and free-exit boundary conditions.

The synthetic scenario uses a non-stationary, divergence-free velocity field to advect smooth initial data over a domain with non-stationary boundary conditions. We use a relatively low-resolution grid here, as we run the global filter for comparison with the distributed filter on this test-case, with the former being more computationally expensive. The grid resolution required for the real data would be too expensive for the global filter.

\subsection{Assimilation of synthetic data}
\label{sec:synthData}

The distributed filter is first tested on the domain, $\Omega = [0,2\pi] \times [0,2\pi]$, discretised into a $10\times10$ element grid with $N=3$. We first generate synthetic observations by solving the advection equation,~\cref{eq:advection}, using the implicit midpoint method with time-step, $\Delta T$, over $t\in[0,T]$ with initial data,
\begin{equation}
\label{eq:synthIC}
c(\bm{x},0)=\sin(x)\cos(y) + 1.2, \quad \bm{x} = (x,y)^\top \in \Omega
\end{equation}
with the divergence-free velocity-field,
\renewcommand\arraystretch{1.2}
\begin{equation}
\label{eq:synthVel}
\bm{u}(\bm{x},t)=\begin{pmatrix}u\\v\end{pmatrix} = \begin{pmatrix}\sin(\frac{x}2)\sin(\frac{y}2)\cos(\frac{2 \pi t}{10}) \\ \cos(\frac{x}2)\cos(\frac{y}2)\cos(\frac{2 \pi t}{10}) \end{pmatrix} 
\end{equation}
and boundary conditions,
\begin{equation}
\label{eq:synthBC}
\begin{split}
c_{lower}(x,t)&=c_{upper}(x,t)= \sin(x)\cos(t), \quad x \in [0,2\pi],\\
c_{left}(y,t)&=c_{right}(y,t)= \sin(y)\cos(t), \quad y \in [0,2\pi].
\end{split}
\end{equation}
The numerical solution, $c(\bm{x}_j,t), \, t \in [\Delta T, T]$ is used as observations, i.e. $y(\bm{x}_j,t)=c(\bm{x}_j,t)+\eta_j(t)$ with $1\%$ signal-to-noise ratio, which is reflected by $r^{tr}=1e-5$ in the matrix $R=r^{tr}I$. We set $\bm{y}_s(t):=\{y(\bm{x}_j,t)\}_{j=1}^{(N+1)^2K}$ and use this vector-function for the tests in the synthetic scenario.

\subsubsection{Global filter}
\label{sec:globalFilter}
For the case of synthetic model parameters and observations, we compare the results of the distributed filter to the global implementation. This involves filtering using a global system matrix as opposed to $K$ local elemental matrices. For this comparison we assume that there is no model error, and boundary conditions are known too, so that $g=0$ and $g_\partial=0$. The global DG system in this case can be written as
\begin{equation}
\label{eq:affineSysGlobal}
\dfrac{d \bm{c}_h}{dt} = A \bm{c}_h + \bm{b}\,, \bm{c}_h(0) = \bm{c}_0+\bm{e}_0^g\,,
\end{equation}
where $A$ is a global system matrix and $\bm{b}$ is a global flux vector. The filter equations are the same as those in~\cref{sec:localMF}, except with global terms in place of local. An apparent advantage of this approach is that unlike in the local formulation, the system matrix, $A$, absorbs the entire inter-element flux term (rhs of~\cref{eq:weakDG}) except in the case where the element boundary in question lies on the domain boundary,~$\Gamma$. In that case, the terms in~\cref{eq:affine} that come from the boundary data contribute to vector, $\bm{b}$. As a result, that flux vector only contains information relating to the domain boundary, while the global matrix, $A$ absorbs everything else. The advantage of this approach is that the global gain, $P$ is informed by richer information than is $P^k$ in the distributed case. However, in the case where observations are incomplete, elements with data neighbouring elements without data may give rise to sharp discontinuities, which will appear in $A$, and then manifest themselves in $P$. In anticipation of this issue, we can define a different global system,
\begin{equation}
\label{eq:affineSysGlobal2}
\dfrac{d \bm{c}_h}{dt} = A_b \bm{c}_h + \bm{b}_b\,, \bm{c}_h(0) = \bm{c}_0+\bm{e}_0^g\,,
\end{equation}
where for the range of rows of $A_b$ relating to a given element, $k$, terms containing the state outside of that element are placed in $\bm{b}_b$. Using this approach, discontinuities due to spatially sparse observations will not appear in $A_b$ and hence neither in the gain, $P$. The gain will now be fed by the same information as in the case of the distributed filter where a similar approach was taken by necessity. The trust in the observations is varied in the same way in the global case as it is for the distributed filter described in~\cref{sec:varyingTrust}.

Later, we compare the cost of the global filter to that of the distributed filter.

\subsubsection{Synthetic Test 1: Full observations over time and space}
\label{sec:fullObsSynth}

In the first test, we initialise the filter to zero, i.e., $\bm{c}_0(\bm{x}) = 0$, and impose the correct boundary conditions,~\cref{eq:synthBC}, and velocity field,~\cref{eq:synthVel}, while taking a single observation over all LGL points at one computational time-level (t=$\Delta T$), where the time step $\Delta T$ is the same for the filter as it is for the forward run for generating observations. We thus set $g=0$ (no model error), $g_0=\frac1{\sqrt{8}}$ (rescaled initial condition error) and $g_\partial=0$ (exact boundary conditions). As a result, $V^k=0$ in~\eqref{eq:ric1}-~\eqref{eq:lMF2}. We also set $q_0=\frac 12$ so that $P^k(0)= (1+2\Delta T)(N+1)^2G_0^k(Q_0^{k})^{-1}G_0^k \approx I$ indicating small trust to the initial conditions. Note, that the minimax filter does not depend on initial conditions, provided it is integrated for a long enough time that $P^k(0)$ is forgotten.

\begin{figure}
        \centering
        \begin{subfigure}[b]{0.49\textwidth}
                \includegraphics[width=\textwidth]{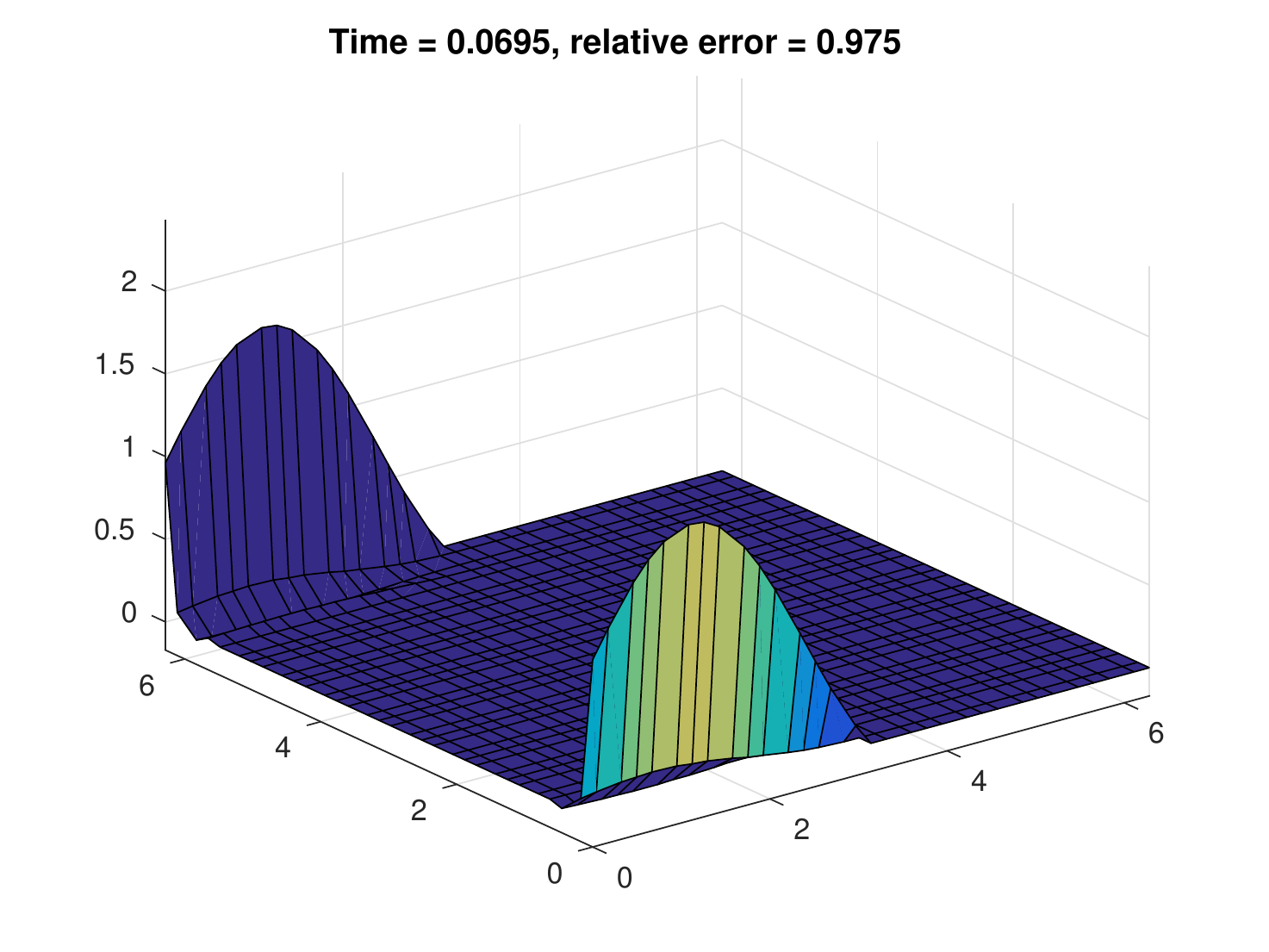}
                \caption{Estimate at $t=\Delta T$: influenced only by initial and boundary conditions  }
                \label{fig:localFilterFullObs1}
        \end{subfigure}%
        ~ %
        \begin{subfigure}[b]{0.49\textwidth}
                \includegraphics[width=\textwidth]{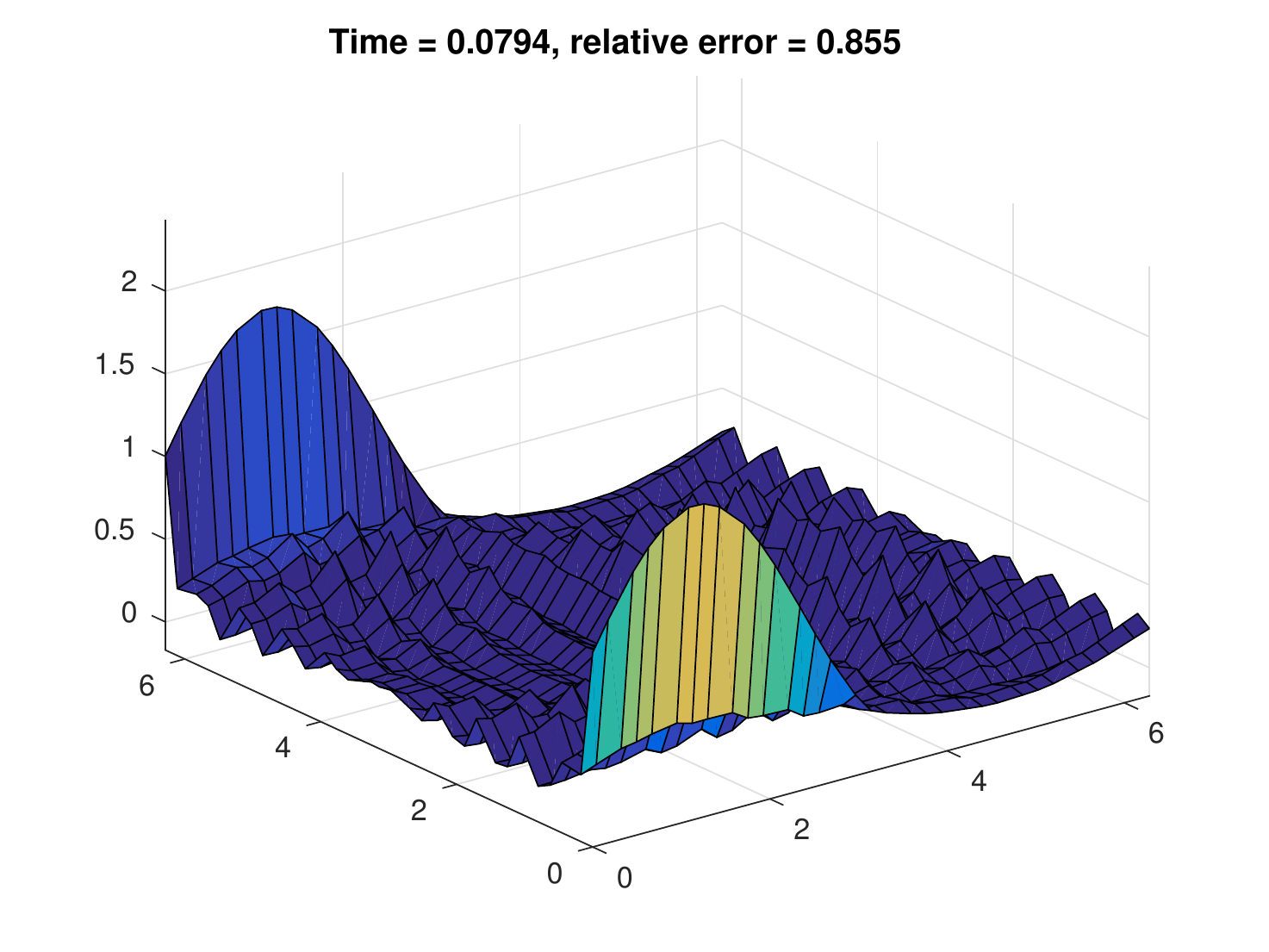}
                \caption{Estimate at $t=0.0794$: assimilating observation}
                \label{fig:localFilterFullObs2}
        \end{subfigure}%
        ~ %

        \begin{subfigure}[b]{0.49\textwidth}
                \includegraphics[width=\textwidth]{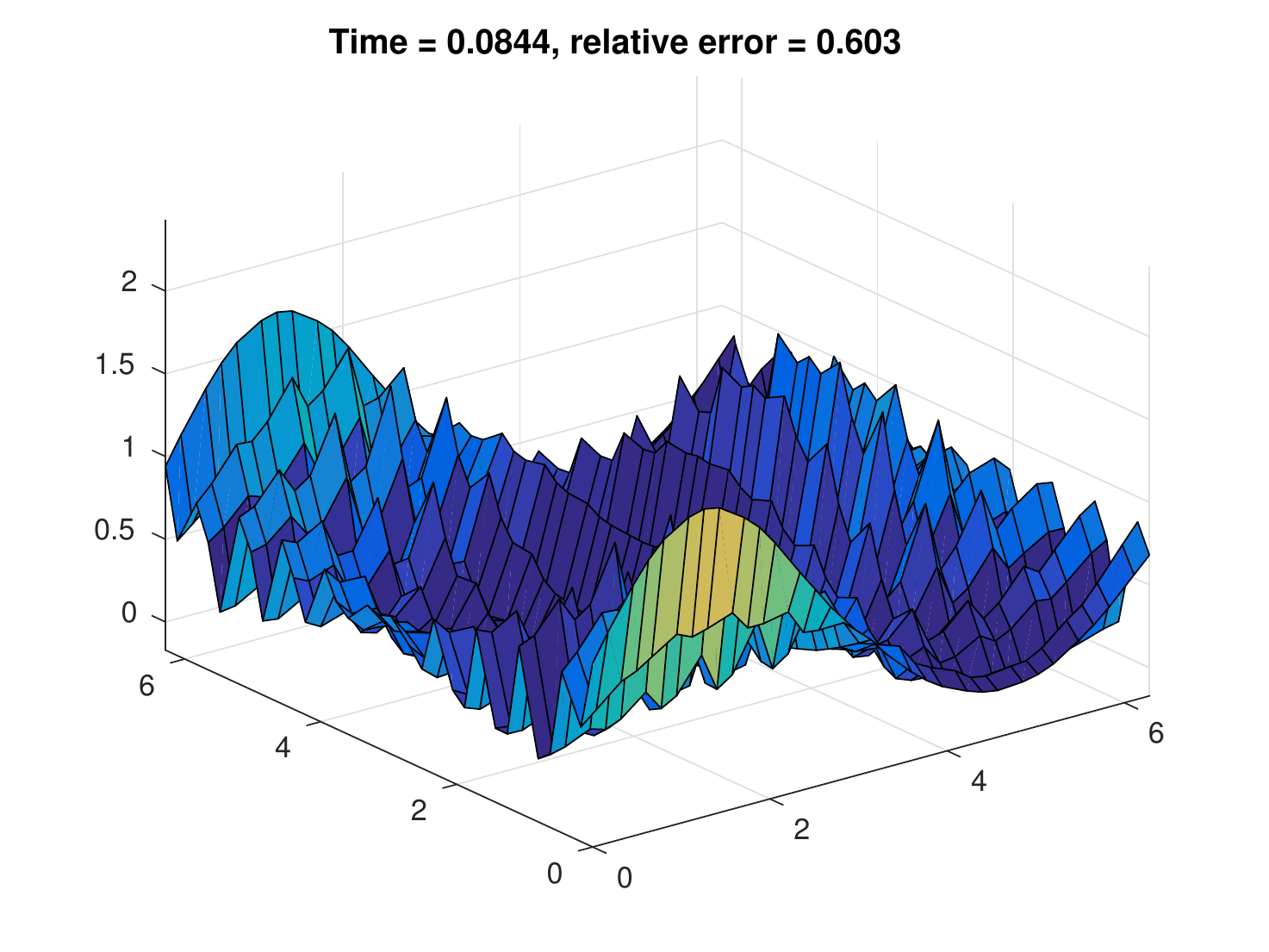}
                \caption{Estimate at $t=0.0844$: assimilating observation}
                \label{fig:localFilterFullObs3}
        \end{subfigure}%
        ~ %
        \begin{subfigure}[b]{0.49\textwidth}
                \includegraphics[width=\textwidth]{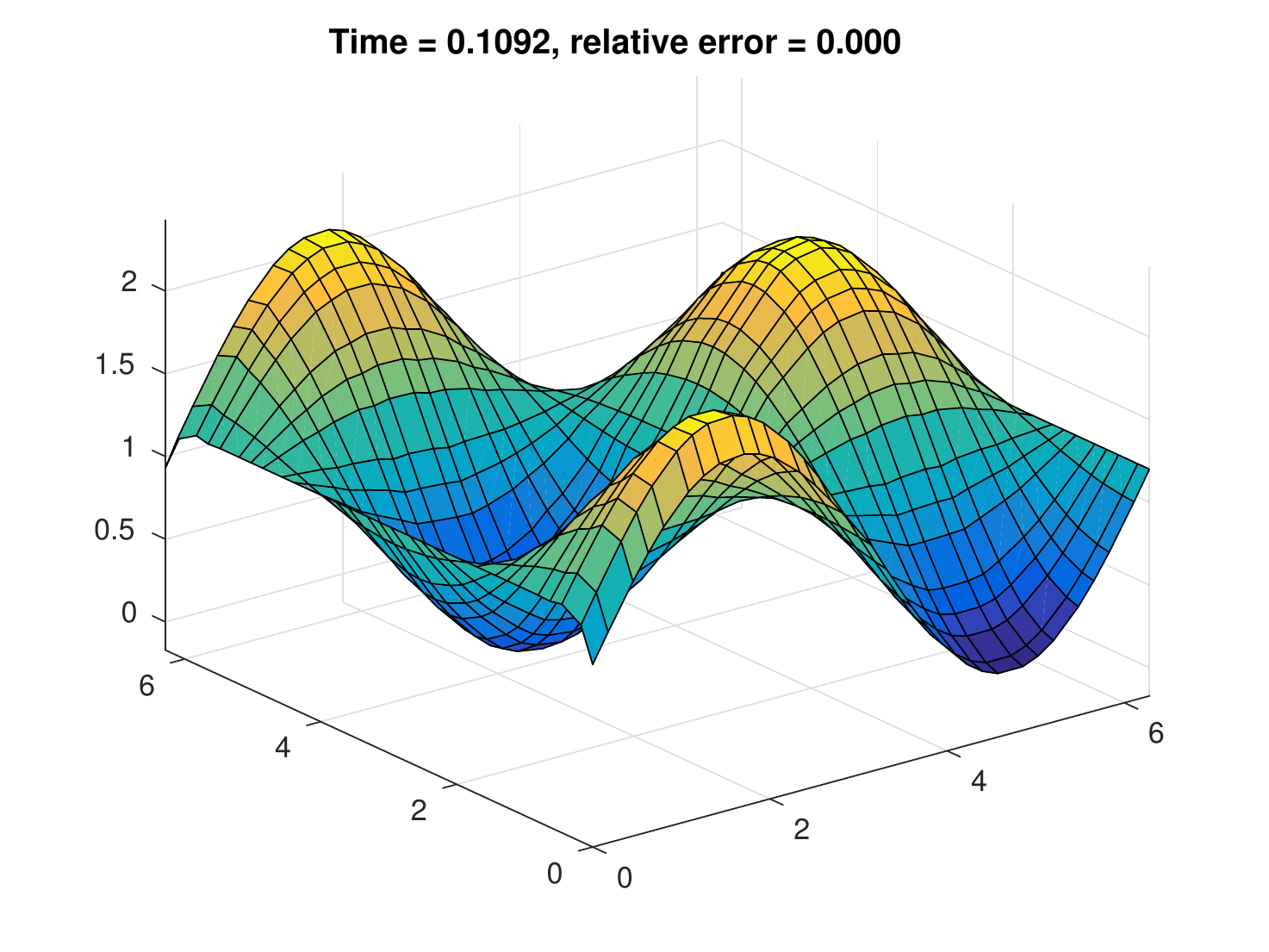}
                \caption{Estimate at $t=0.1092$: assimilated observation}
                \label{fig:localFilterFullObs4}
        \end{subfigure}%

    \caption{Assimilation of a single full observation with precise boundary conditions and velocity field (Synthetic Test 1)}\label{fig:localFilterFullObs}
\end{figure}

\begin{figure}[htbp]
  \centering
  \label{fig:a}\includegraphics[width=0.6\textwidth]{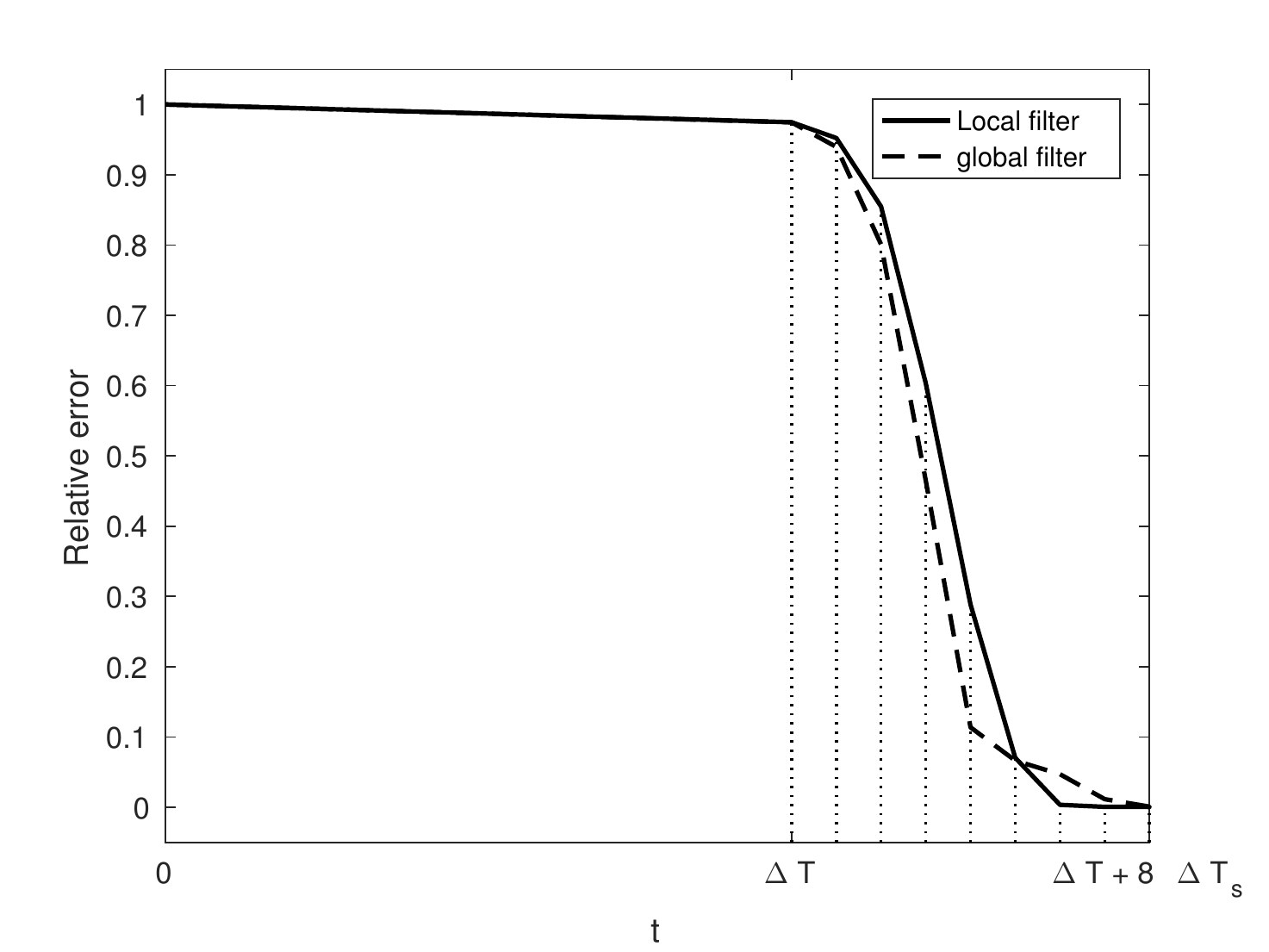}
  \caption{Estimate error for local and global filters for 1 full observation with precise boundary conditions and velocity field (Synthetic Test 1)}
  \label{fig:fullObsCorrectBCErrors}
\end{figure}

This is a simple experiment to observe how the filter reacts to the incorrect initial condition. The time-evolution of the distributed filter for for $t\in[0,0.1092]$ is shown in figure~\cref{fig:localFilterFullObs}.
The error shown in the figures is the relative $L^2$ error, i.e. $\| \bm{y}(t) - \hat{\bm{c}}_h(t) \| / \| \bm{y}(t) \|$. The computational time-step, $\Delta T = 0.0695$ and the filtering time-step, $\Delta T_s = \Delta T/14$ meaning we run the filter with varying trust 14 times over $t\in[\Delta T, 2\Delta T]$ for the observation, $\bm{y}_s(\Delta T)$, which becomes available at $t=\Delta T$ (see~\cref{sec:varyingTrust}~and~\cref{alg:localFilter}). Until $t=\Delta T$, the estimate is zero everywhere apart from at the regions affected by the boundary conditions, which may induce an inflow depending on the direction of the velocity field on the boundaries (see~\cref{sec:DGform}). We see this in~\cref{fig:localFilterFullObs1}, where observations have not yet become available. In~\cref{fig:localFilterFullObs2}, the first observation is being assimilated at time $t = \Delta T + 2\Delta T_s$, i.e., the 2nd time-step of the inner time-loop in~\cref{alg:localFilter}. By $t=0.1092$, the first observation has been assimilated and the relative error is~$\approx 1e-4$. This is then compared to the result of the global filter with fully global system matrix, $A$, and boundary flux vector, $\bm{b}$ (see~\cref{eq:affineSysGlobal}), as described in~\cref{sec:globalFilter}, under the same conditions. The relative errors for both the local and global filters over time are shown in~\cref{fig:fullObsCorrectBCErrors}, where the observation becomes available at $t=\Delta T$. 
At $t=0.1092$ (i.e, $t=\Delta T + 8\Delta T_s$), the relative error in the global case is also~$\approx 1e-4$.

\subsubsection{Synthetic Test 2: Sparse observations in space}
\label{sec:sparseObsSynth}

\begin{figure}
        \centering
        \begin{subfigure}[b]{0.49\textwidth}
                \includegraphics[width=\textwidth]{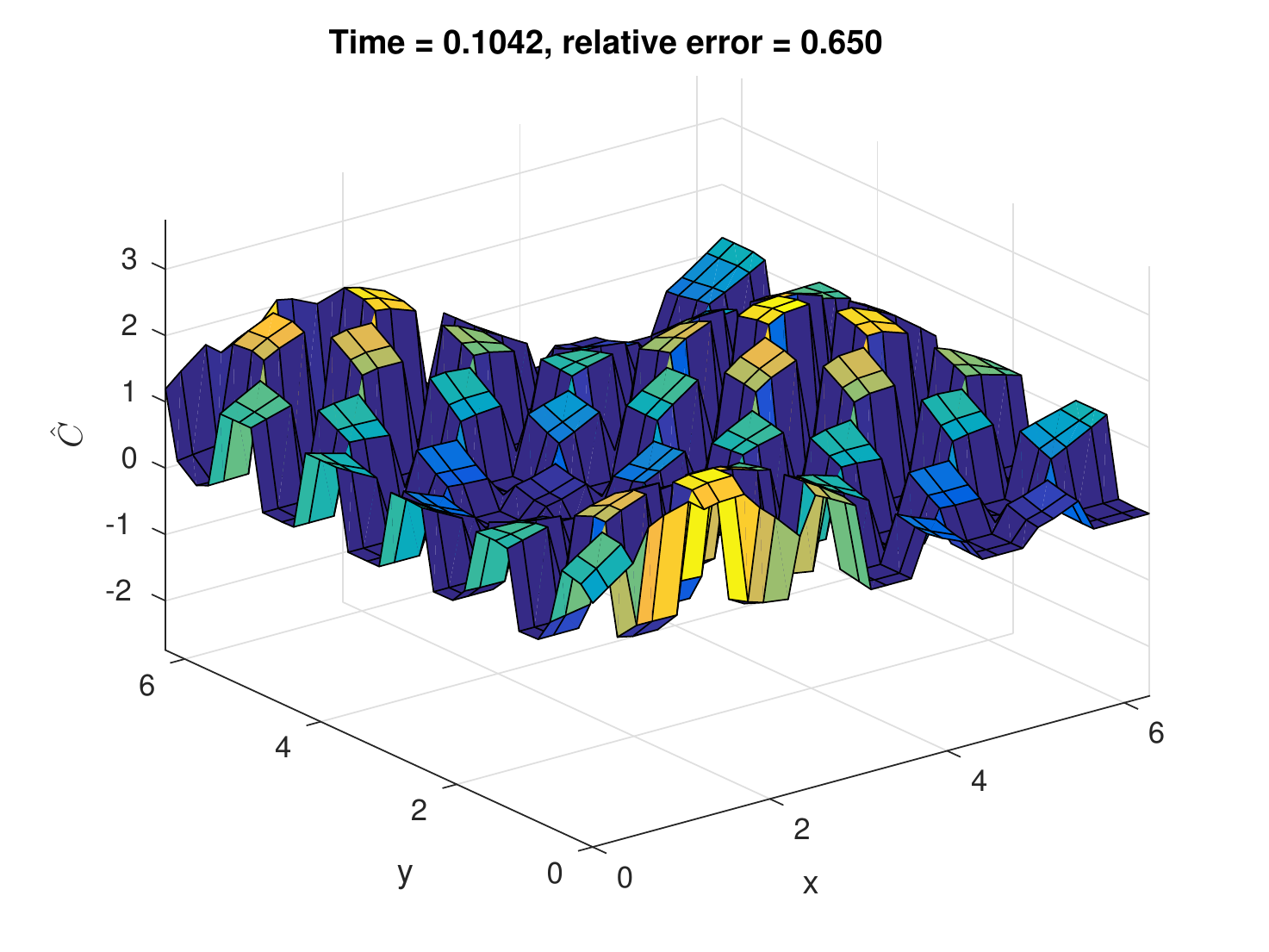}
                \caption{Estimate at $t=0.1042$}
                \label{fig:localFilterPartialObs1}
        \end{subfigure}%
        ~ 
        \begin{subfigure}[b]{0.49\textwidth}
                \includegraphics[width=\textwidth]{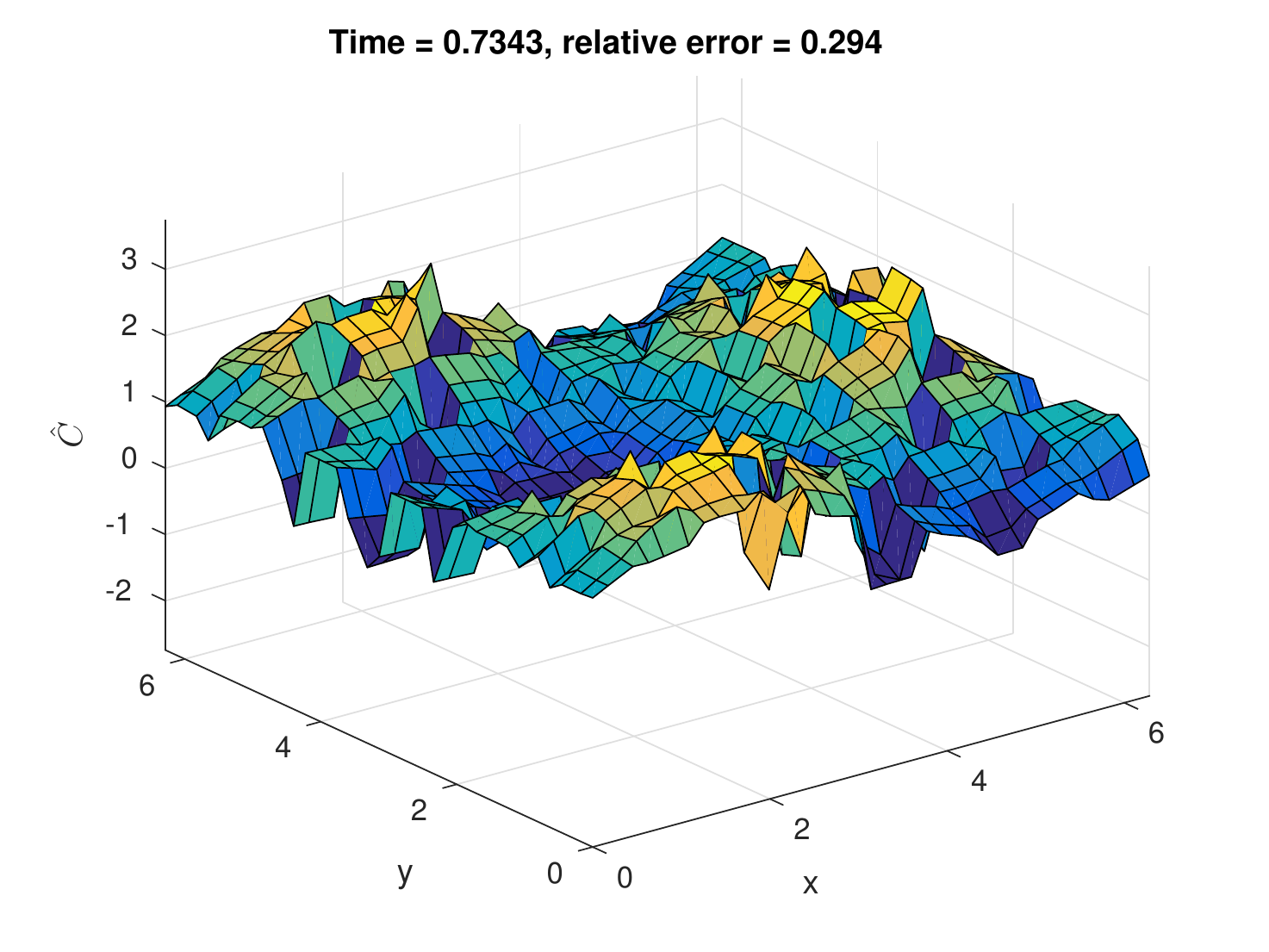}
                \caption{Estimate at $t=0.7343$}
                \label{fig:localFilterPartialObs2}
        \end{subfigure}%

        \begin{subfigure}[b]{0.49\textwidth}
                \includegraphics[width=\textwidth]{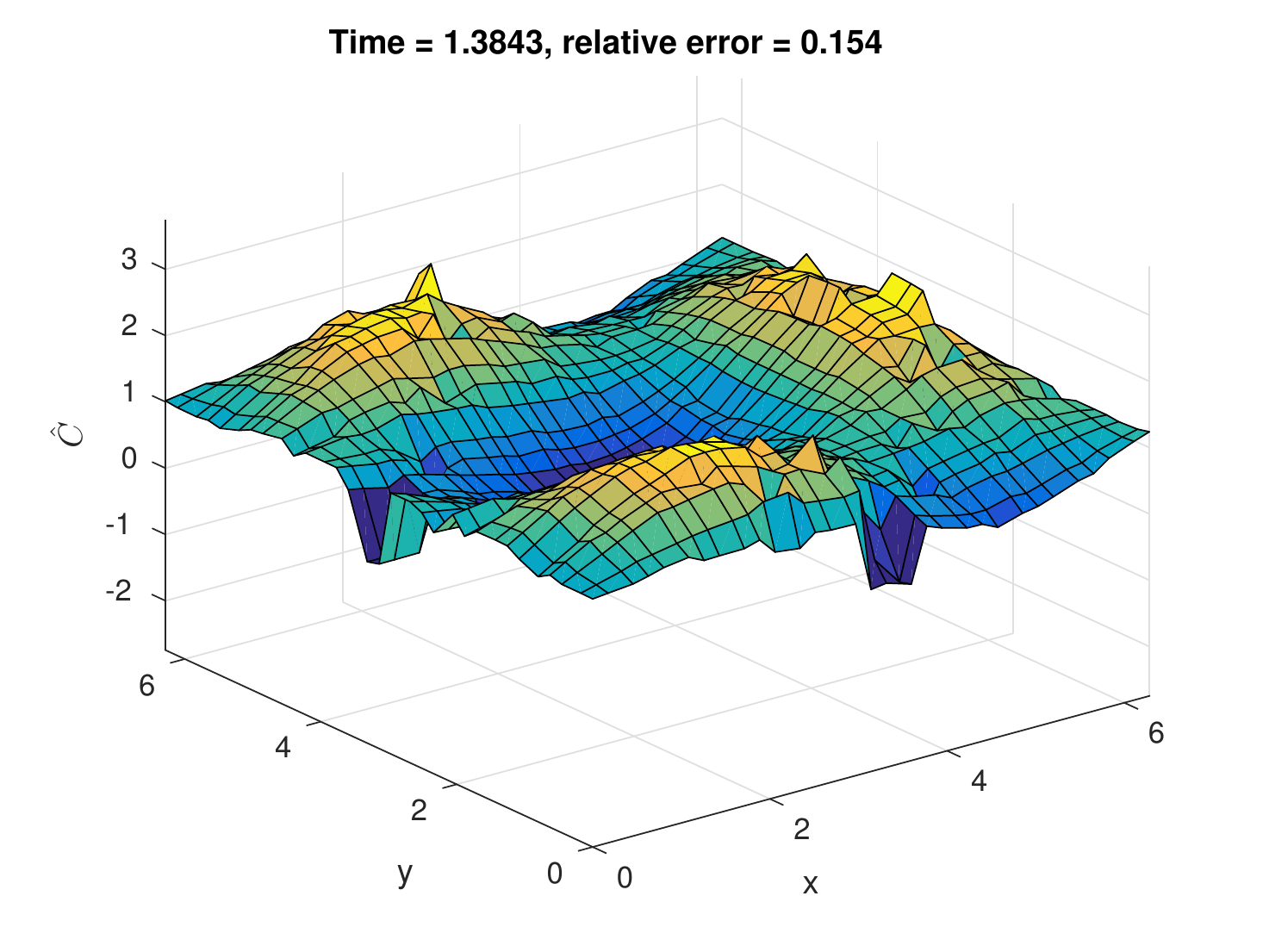}
                \caption{Estimate at $t=1.3843$}
                \label{fig:localFilterPartialObs3}
        \end{subfigure}%
        ~ 
        \begin{subfigure}[b]{0.49\textwidth}
                \includegraphics[width=\textwidth]{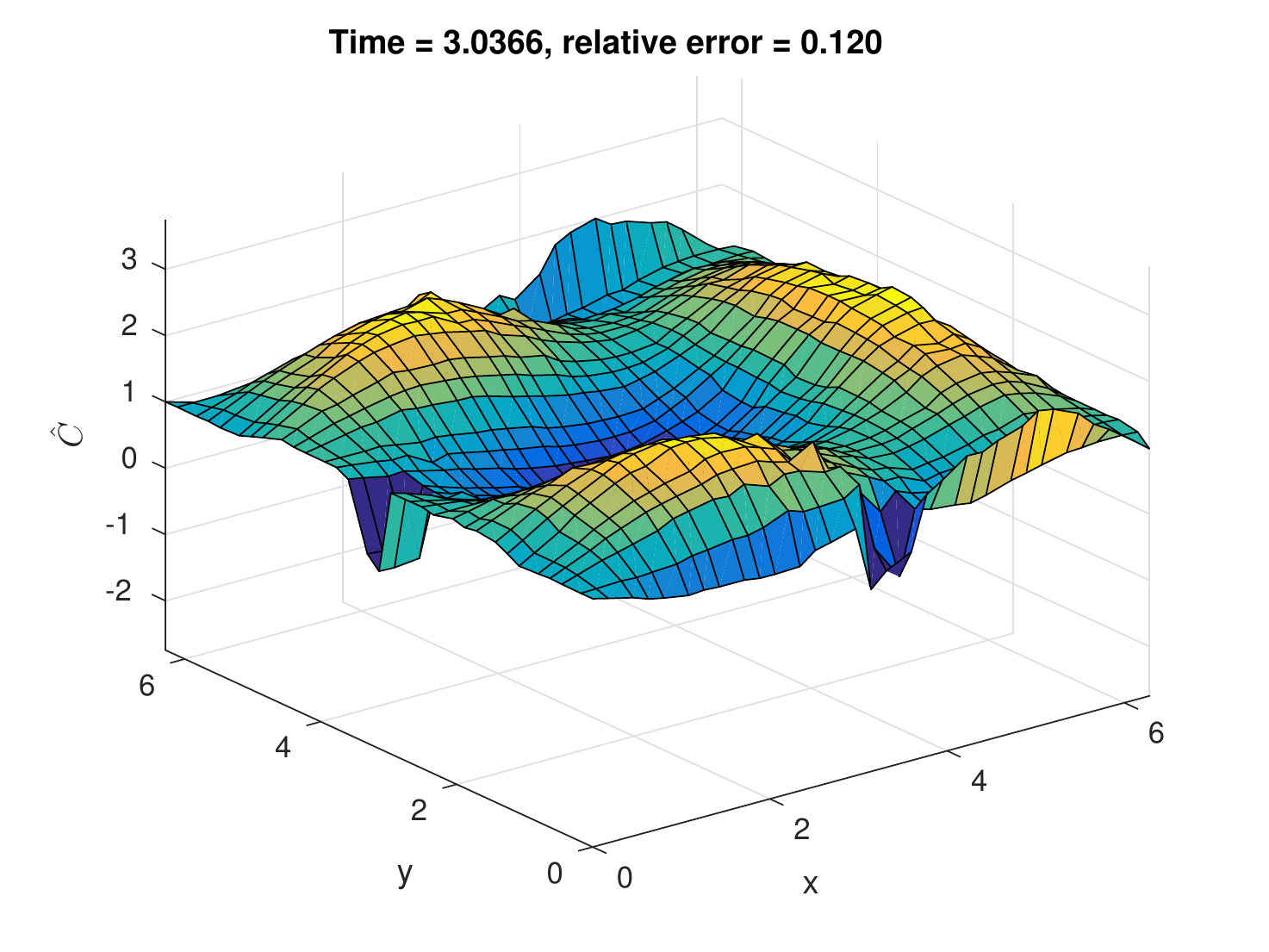}
                \caption{Estimate at $t=3.0366$}
                \label{fig:localFilterPartialObs4}
        \end{subfigure}%

    \caption{Assimilation of a partial observations with precise boundary conditions and velocity field (Synthetic Test 2)}\label{fig:localFilterPartialObs}
\end{figure}

\begin{figure}
        \centering
        \begin{subfigure}[b]{0.49\textwidth}
                \includegraphics[width=\textwidth]{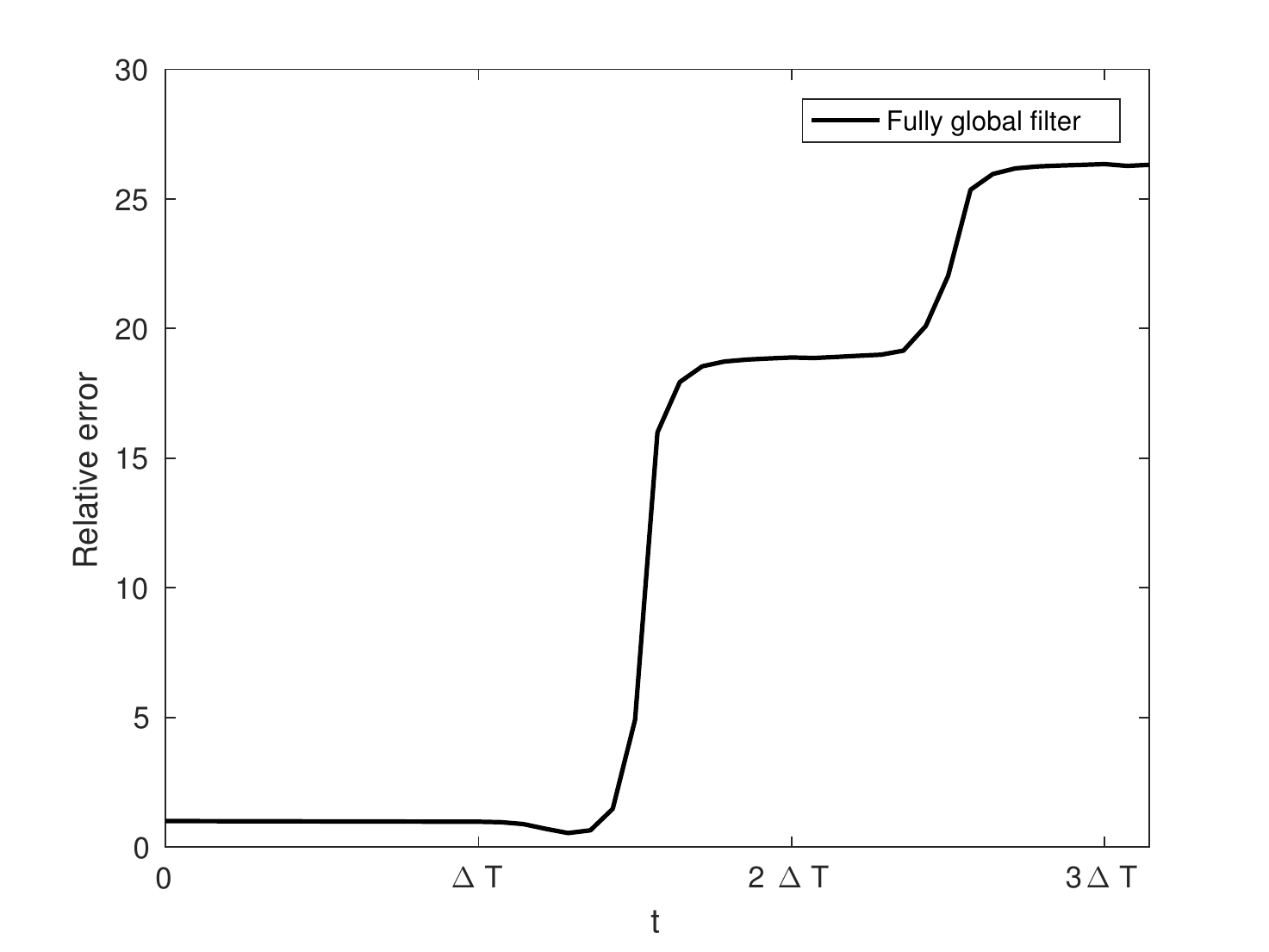}
                \caption{Relative error for fully global filter \newline}
                \label{fig:partialObsCorrectBCErrorGlobal}
        \end{subfigure}%
        ~ 
        \begin{subfigure}[b]{0.49\textwidth}
                \includegraphics[width=\textwidth]{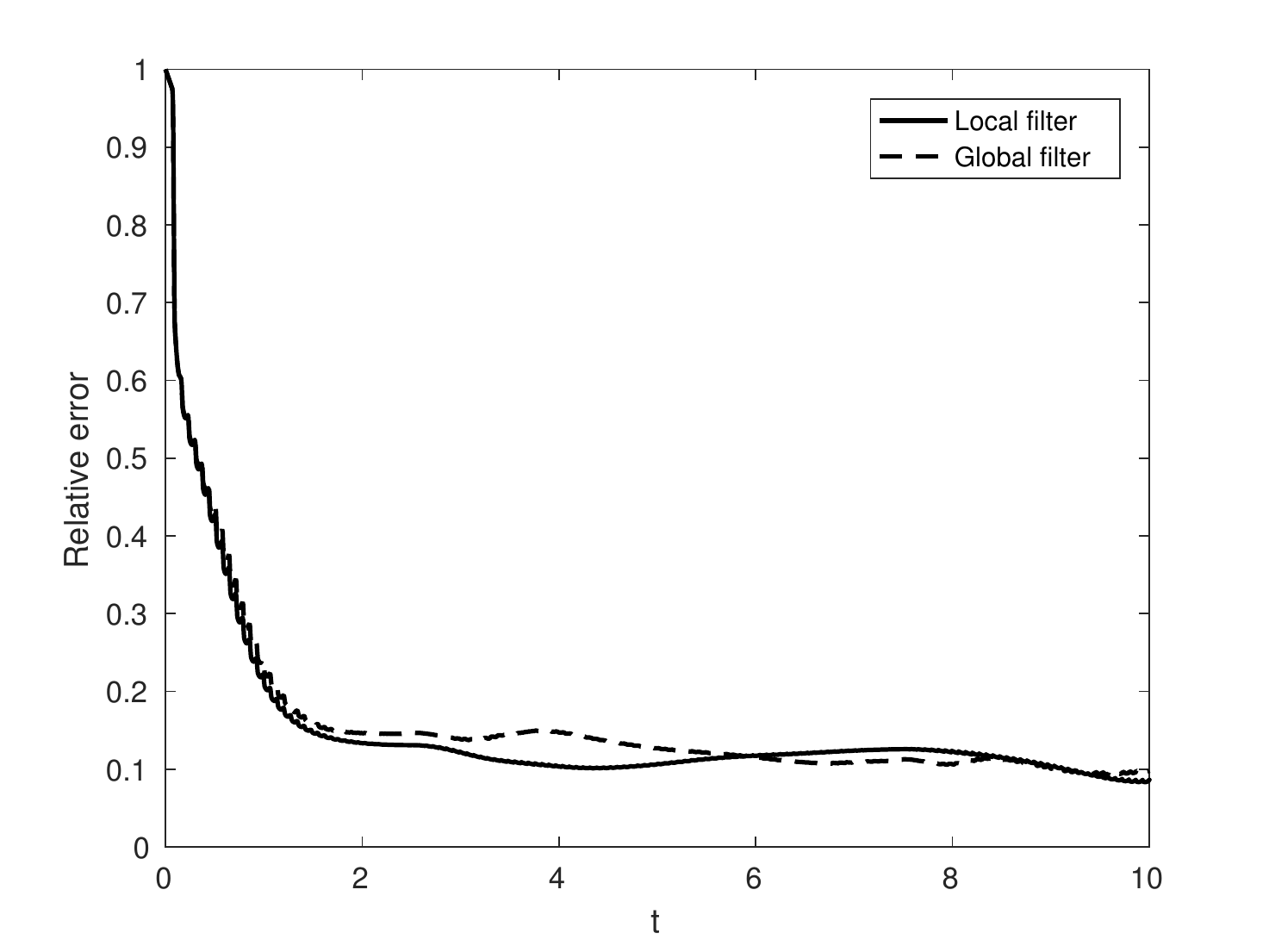}
                \caption{Relative error for local and modified global filters}
                \label{fig:partialObsCorrectBCErrorGlobalAndLocal}
        \end{subfigure}%
    \caption{Relative error for partial observations with precise boundary conditions and velocity field  (Synthetic Test 2)}\label{fig:partialObsCorrectBCErrors}
\end{figure}

\begin{figure}[htbp]
  \centering
  \includegraphics[width=0.6\textwidth]{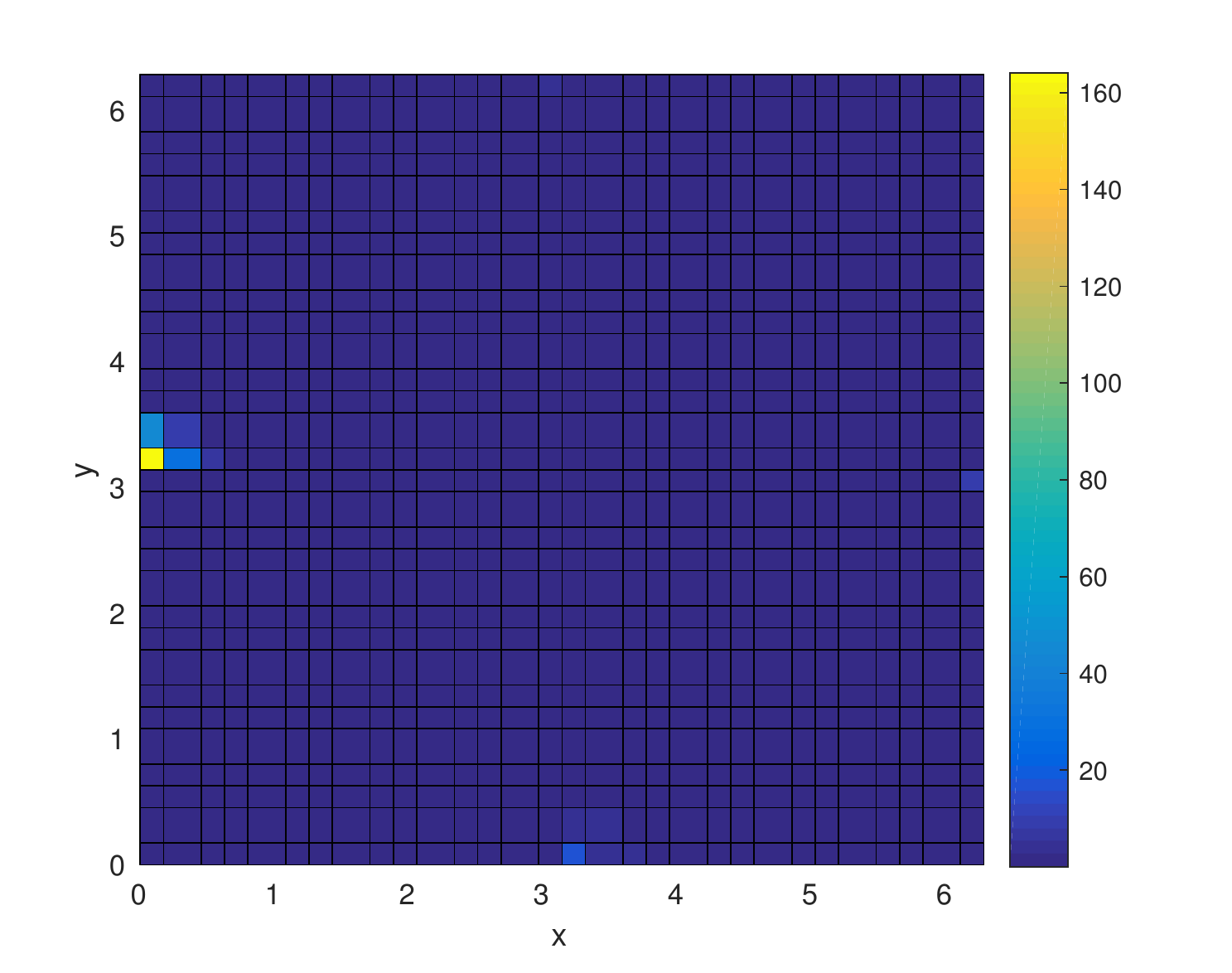}
  \caption{Worst-case estimation error for partial observations with precise boundary conditions and velocity field at $t=10$ for distributed/global filters (Synthetic Test 2)}
  \label{fig:partialObsCorrectBCWCError}
\end{figure}

In the second test, we provide the filter with the correct boundary conditions,~\cref{eq:synthBC}, and velocity field,~\cref{eq:synthVel}, but incomplete observations. Specifically, we equip every other element with observations in a ``chequered'' pattern. This experiment is designed to test the filter in the presence of discontinuities in the observations that may occur when dealing with sparse data. Like in the previous experiment, the estimate is initialised to zero, while observations are available at times, $t=\Delta T, 2\Delta T, \ldots$. We use the same $g$, $g_0$, $g_\partial$ and $q_0$. The time-evolution of the estimate is shown in~\cref{fig:localFilterPartialObs}.
We see the first observation being assimilated in~\cref{fig:localFilterPartialObs1} where the ``chequered'' observation pattern is apparent. The subsequent figures show the estimate as further observations are assimilated; we see that over time, the observation pattern becomes less apparent as the velocity field induces a flux between neighbouring elements, and the relative error decreases over time. We repeat the same test with the global filter on the system,~\cref{eq:affineSysGlobal}, for comparison. This fails to converge, as we see from~\cref{fig:partialObsCorrectBCErrorGlobal}. 
As discussed in~\cref{sec:globalFilter}, discontinuities generated by sparse observations will manifest themselves in the global gain, possibly leading to instabilities. This could account for the failure of the global filter in this case. We re-attempt the experiment with the system,~\cref{eq:affineSysGlobal2}, where the system matrix, $A_b$, and flux vector, $\bm{b}_b$, are described in~\cref{sec:globalFilter}. To summarise that discussion, the system matrix is constructed as follows: for a given element, components of the state on neighbouring elements are placed in the flux vector, $\bm{b}_b$ instead of in $A_b$. As a result, the latter will not ``see'' the discontinuities induced by sparse observations, and consequently, neither will the global gain, $P$. In this case, the performance of the global filter is similar to that of the distributed filter, as seen in~\cref{fig:partialObsCorrectBCErrorGlobalAndLocal}, where, at $t=10$, the relative $L^2$ error is $\approx 0.09$ for the distributed filter and $\approx 0.1$ for the global filter. This is not surprising, since the gain in either case is influenced by the same information, as explained in~\cref{sec:globalFilter}. An important point to make is that the apparent advantage of the fully global filter (discussed in~\cref{sec:globalFilter}) is gone when discontinuities due to observation sparsity are present. While having the system matrix absorb all but the boundary terms is generally preferable because of the richer information available to the filter, we saw that in the case of sparse observations, the discontinuities generated may cause the filter to fail. We also note that in the experiment with full observations, the fully global filter performed well, which supports our assertion that the presence of discontinuities due to observation sparsity was responsible for the failure of that filter. In~\cref{fig:partialObsCorrectBCWCError}, we see the worst-case error for the full domain at $t=10$ for the local filter. The relative $L^2$ error between this and the worst-case error for the global filter were found to be within $1e-8$ of one another. The higher worst-case errors at some points on the boundaries are due to the advection field being close to zero at those points. This results in little transport in those regions, which consequently results in little element-to-element communication, and thus higher errors. We conclude that for this specific case, the proposed distributed filtering strategy yields a similar performance to that of the global filter while running at a much lower computational cost.    

\subsubsection{Synthetic Test 3: Sparse observations in space with imprecise knowledge of boundary conditions and velocity field}
\label{sec:sparseObsSynthWrongBC}

\begin{figure}
        \centering
        \begin{subfigure}[b]{0.49\textwidth}
                \includegraphics[width=\textwidth]{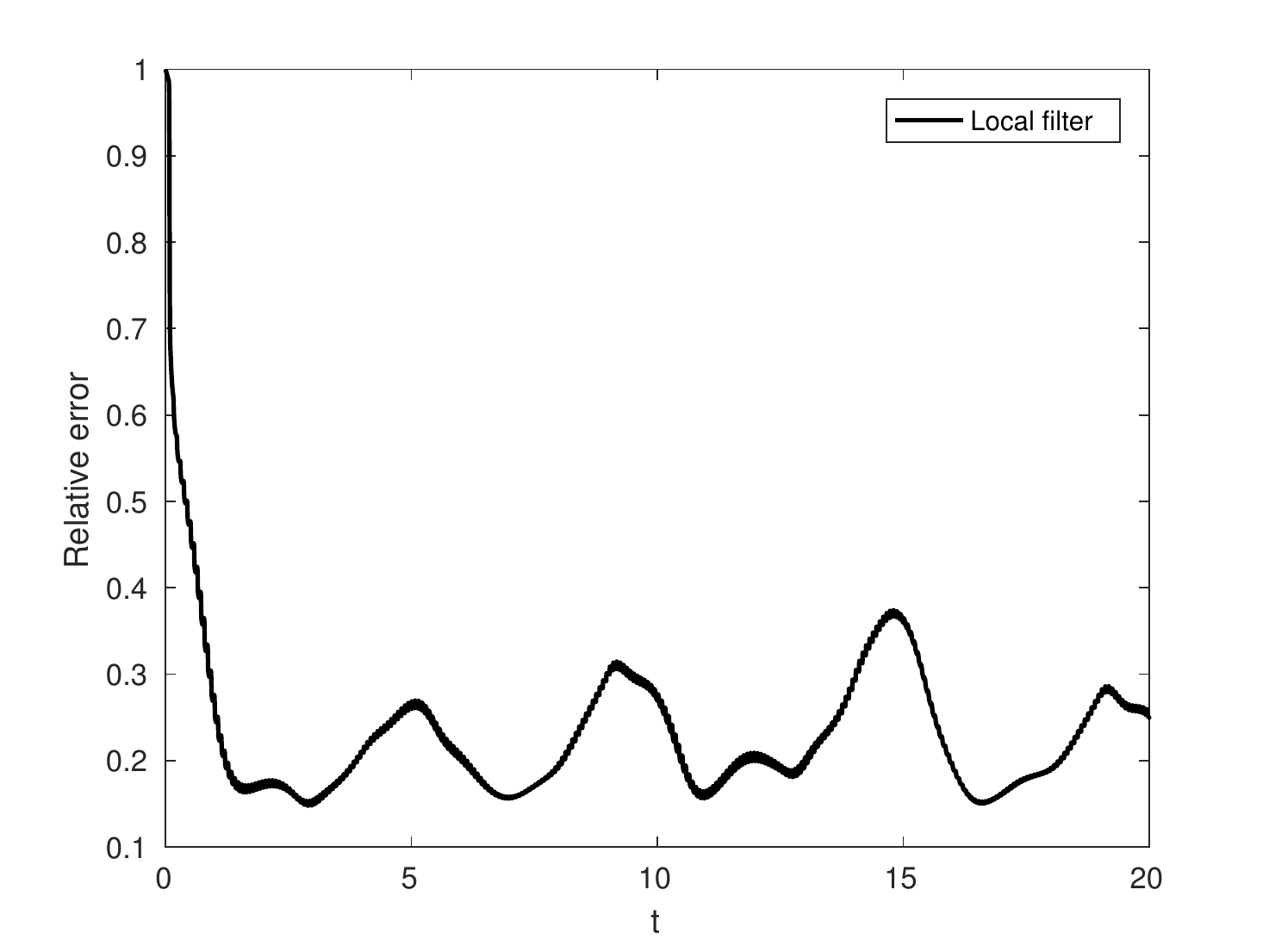}
                \caption{Relative error for distributed filter \newline}
                \label{fig:partialObsIncorrectBCErrors1}
        \end{subfigure}%
        ~
        \begin{subfigure}[b]{0.49\textwidth}
                \includegraphics[width=\textwidth]{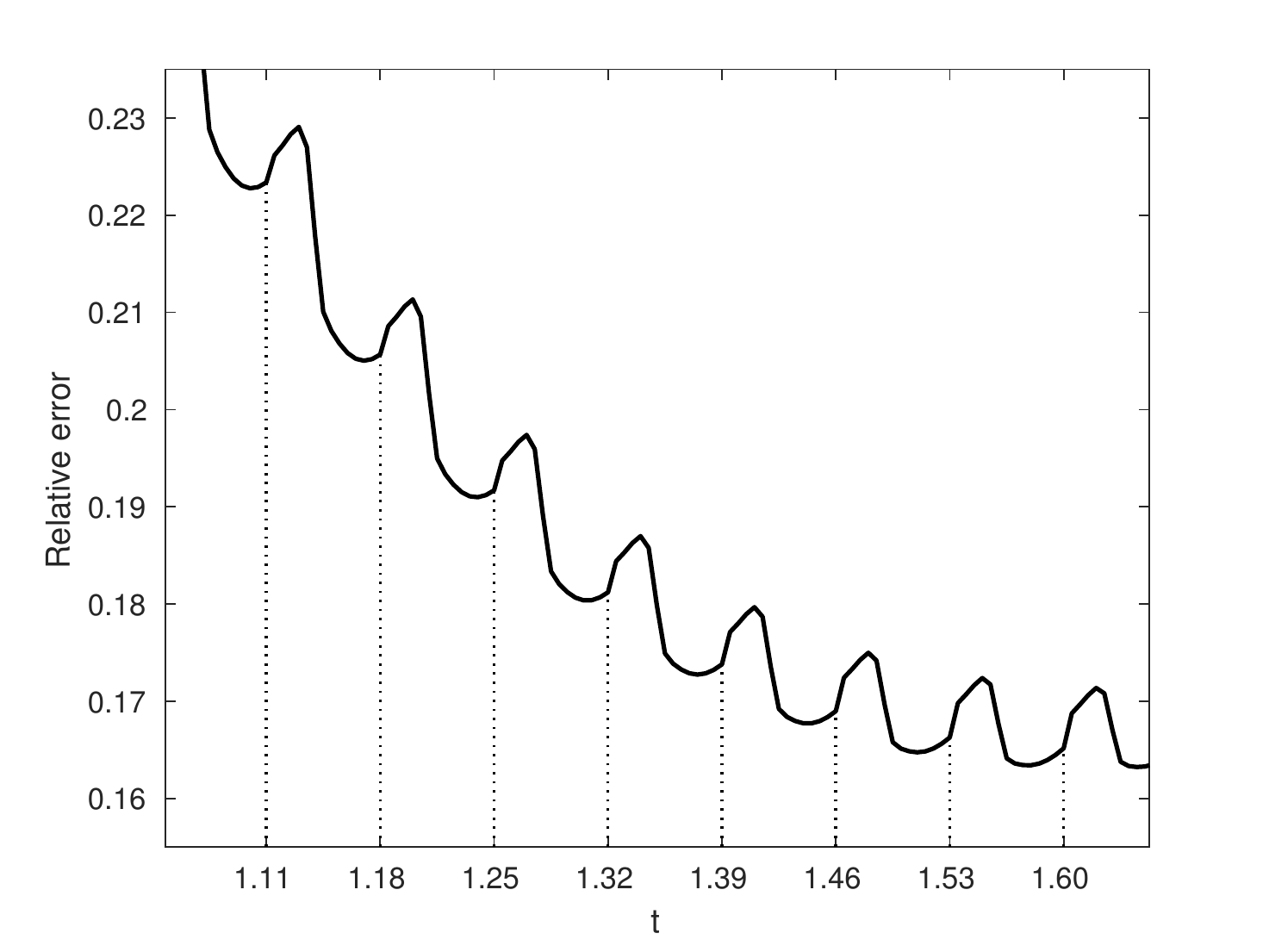}
                \caption{Relative error for distributed filter over short time}
                \label{fig:partialObsIncorrectBCErrors2}
        \end{subfigure}%
    \caption{Relative error for partial observations with imprecise boundary conditions and velocity field  (Synthetic Test 3)}\label{fig:partialObsIncorrectBCErrors}
\end{figure}

\begin{figure}[htbp]
  \centering
  \includegraphics[width=0.6\textwidth]{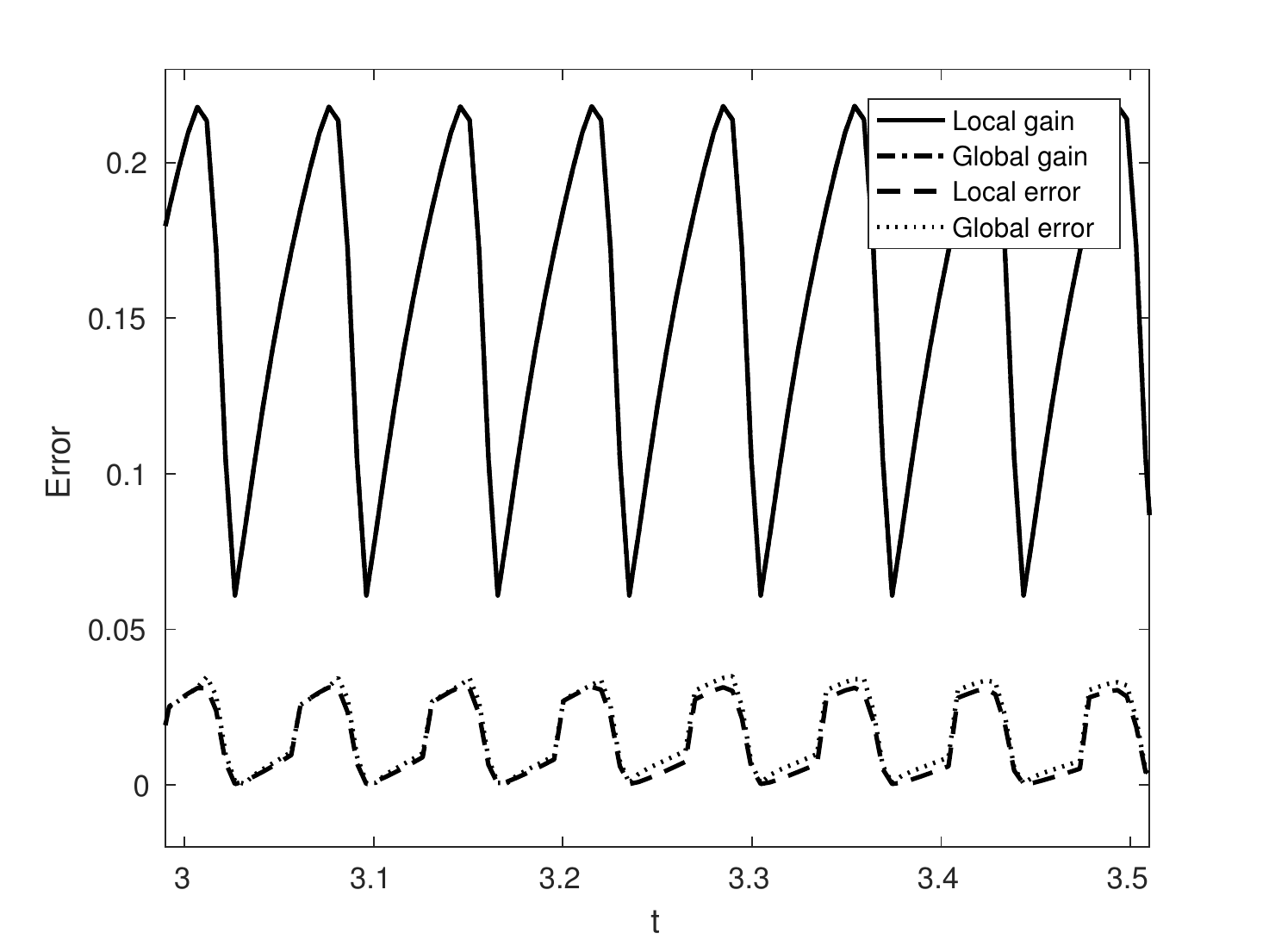}
  \caption{Square root of the worst-case estimation error and estimate error at single node for partial observations with imprecise boundary conditions and velocity field (Synthetic Test 3)}
  \label{fig:partialObsWrongBCGainAndError}
\end{figure}

The next synthetic experiment we perform is with the same ``chequered'' observations but with imprecise boundary conditions and advection velocity field. To do this, we use~\cref{eq:synthBC} and~\cref{eq:synthVel} with a time-shift, $t_s$. In this way, the advection field and boundary conditions for the filter are out of phase with those used to generate the observations, $\bm{y}_s$. We set $g=g_\partial=1$ and $q=q_\partial=(N+1)^2$. In~\cref{fig:partialObsIncorrectBCErrors1}, we see the relative error of the distributed filter over time when we run the filter with a time-shift of $t_s=-1.5$ for both the advection field and boundary conditions.
Not only is the error generally higher than in the case of \emph{Synthetic Test 2}, it is also seen to fluctuate over time. This appears to be due to the boundary conditions and advection field, which are periodic in time. This emphasises the impact of incorrect model parameters on the estimate; although high trust is placed in the observations, the incorrect model steers the estimate off track. The filter appears to perform better when the difference between the assumed and real boundary conditions and advection field are small, which in this example occurs periodically. In~\cref{fig:partialObsIncorrectBCErrors2}, the error over a relatively short time interval is shown, where the times at which data become available are indicated by vertical dashed lines. Here, we can see the change in error as the trust is varied over each time-step, and how this is counteracted by the model, which steers the estimate away from the truth. The absolute error at a single node for both the local and global filters is shown in~\cref{fig:partialObsWrongBCGainAndError} for a short time interval. Also shown is the square root of the worst-case estimation error~\eqref{eq:lMF_err}, i.e. the square root of the diagonal entry of the gain $P^k$ for the node in question.
We see that both the gains and the errors for local and global are close, as expected, and that the error is bounded by the square root of the gain entry as required.

\subsubsection{Computational complexity and scalability}
\label{sec:highResTest}

Here we provide a high level description of the computational cost measuring the latter in terms of the number of linear solves required to make an estimate for one time step. In the case of the global filter we need to perform one sparse linear solve with a matrix of dimension $2(N+1)^2*K$, where $N$ is the degree of the Lagrange polynomials, and $K$ is the number of elements used for the DG method. The matrices $A$ and $A_b$ are however, sparse, with the number of non-zero terms~$<1\%$, so a fast algorithm for sparse matrices such as GMRES~\cite{saad1986gmres} could be used for linear solves involving this matrix. However, the matrices $U$ and $V$ in~\cref{eq:UVj} are dense, and the computation of $P=VU^{-1}$ at each time-level requires $K(N+1)^2$ linear solves with the matrix $U$, which is of dimension $K(N+1)^2 \times K(N+1)^2$. So, effectively, we need to invert $U$ which costs at least $O((N+1)^6 K^3)$ operations (for Gaussian elimination type linear solvers). Finally, to obtain the estimate using~\cref{eq:filter}, a matrix of dimension $K(N+1)^2\times K(N+1)^2$ containing the dense matrix $P$, must be inverted. In contrast, to compute the distributed filter at each element, one needs to perform the same number of linear solves as the global filter but with matrices whose dimension is a factor of $K$ times smaller: namely, one sparse linear solve with a matrix of dimension $2(N+1)^2\times 2(N+1)^2$ and so on. In particular, in this case inverting $U$ costs $O((N+1)^6)$ so, in total, we need $O((N+1)^6K)$ inversions. As a result, the total cost of computing the distributed filter scales linearly with $K$ whereas for the global case the total computational cost is at least polynomial in $K$. Note that it is also possible to devise a hybrid or ``semi-local'' approach whereby neighbouring elements are grouped into regions on which a system matrix $A^k$ is assembled from the elemental mass and stiffness matrices of the elements in the regions. In this case the computational cost of the distributed filter will be different.

\begin{figure}[htbp]
  \centering
  \includegraphics[width=0.8\textwidth]{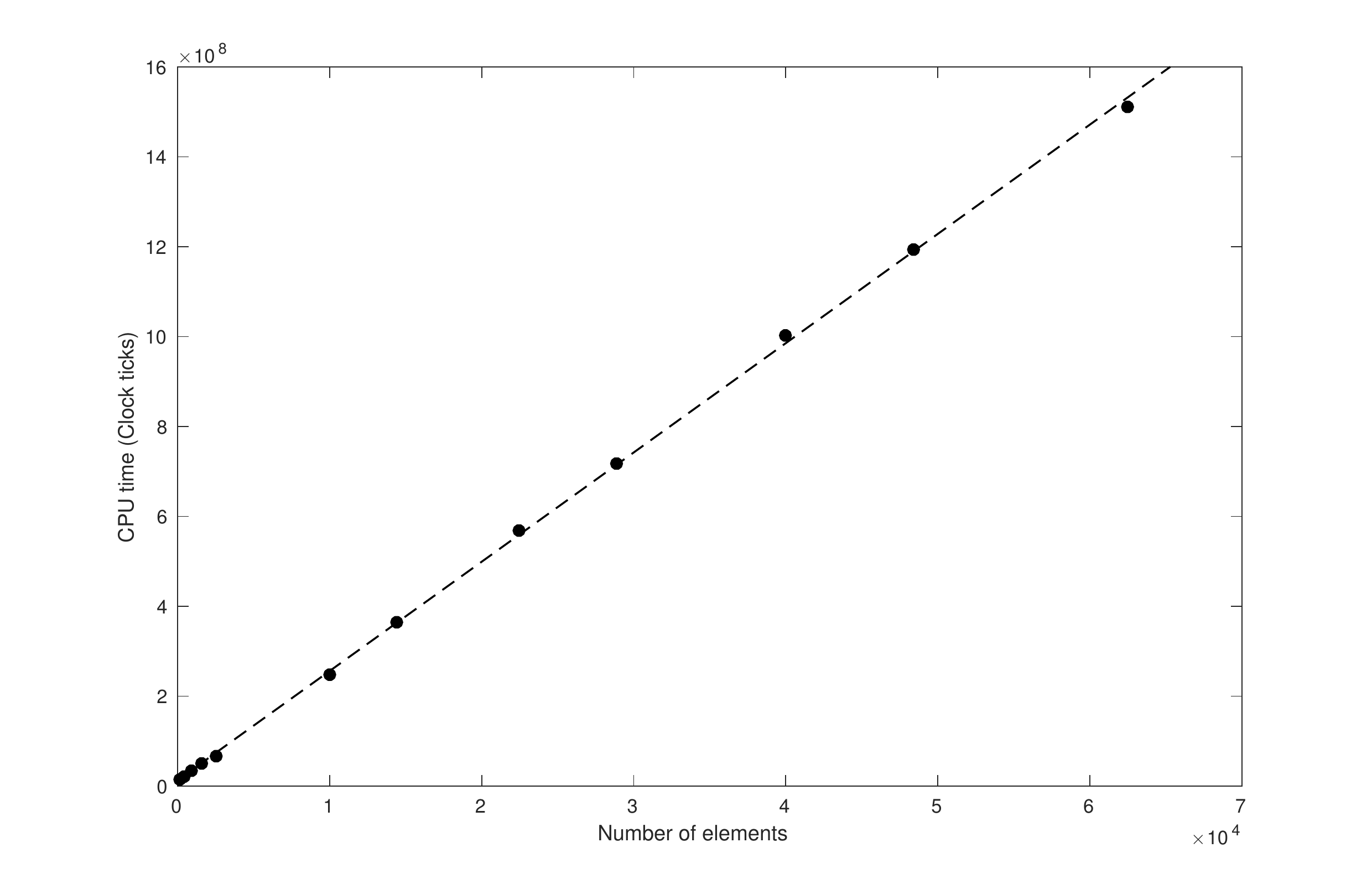}
  \caption{CPU time taken to assimilate one observation vs number of elements in DG grid}
  \label{fig:timings}
\end{figure}

To illustrate the scalability of the algorithm, we measure the CPU time taken to carry out the assimilation of a single set of observations for grids of different resolutions. This process requires the filtering to be performed $n_s$ times in order to ramp up the trust in the observations (see~\cref{sec:varyingTrust}). In these experiments, $n_s=14$. The results are shown in \cref{fig:timings}. 
We see that the cost scales linearly with the number of elements, as expected. The highest resolution grid depicted in the figure has $250 \times 250$ elements with $N=3$ (i.e., state of dimension $10^6$). For that grid, the experimental setup was the same as that in `Synthetic Test 2' described in \cref{sec:sparseObsSynth} (i.e., $10 \times 10$ grid of alternating observation regions and knowledge of advection field and boundary conditions). The filter was run until $t=1.5$ by which time the relative error was $\approx 0.17$. This demonstrates the efficacy and tractability of the algorithm at high resolution. We implemented~\cref{alg:localFilter} in \verb!C++!, and parallelised the element loop (over $k$ in~\cref{alg:localFilter}) using OpenMP. The numerical experiments have been conducted on an IBM POWER8 machine with 196 cores and 0.5TB of RAM.

\subsection{Real data scenario: satellite image assimilation}
\label{sec:COD}

The distributed filtering algorithm is next tested on real data. For this experiment, the observations consist of satellite images depicting cloud optical depth, available at 15-minute intervals. We use a domain that spans 16 degrees of longitude and 12 degrees of latitude which, in the region in question, covers approximately $1.68e6$ square kilometres. The domain is discretised into a $70 \times 70$ element grid with $N=3$, which is close to the resolution of the images. The resulting system is of dimension $19600$, which is considerably higher than that in the synthetic case where it was $400$. 

\subsubsection{Satellite image data and advection field}
\label{sec:CODdata}

Using a purely synthetic velocity field to advect the satellite images is not desirable as the resulting motion may not appear natural. Instead, we compute a sequence of fields that captures the motion of the sequence of images approximately by employing an optical flow estimation procedure~\cite{optflow_IJCV14}. We will not describe the basic optical flow estimation here but note that the procedure does not necessarily produce a divergence-free velocity field. In order to obtain a divergence-free field, we perform a projection and reconstruction of the velocity using the standard vorticity-stream formulation commonly used to solve the incompressible (divergence-free) Navier Stokes equation. This procedure is described briefly below. 

The scalar field, vorticity, $\xi=(\nabla \times \bm{u})\cdot\bm{e}_z$, is obtained numerically from the optical velocity field and is then projected into the space, $\mbox{span}\{\phi_n(x)\phi_m(y)\}_{|n|,|m| \leq N_f/2}$ of complex exponentials, yielding $N_f+1$ projection coefficients, $a_{nm}$, where $|n|,|m| \leq N_f/2$. Defining the stream function, $\Psi$ using the Poisson equation, $-\Delta \Psi = \xi$, we can express the velocity components as $u=\Psi_y$ and $v=-\Psi_x$. Decomposing the Laplacian operator by taking complex exponentials as eigenfunctions, we can reconstruct the velocity field using the corresponding eigenvalues, $\lambda_{nm}$, as follows: 
\begin{equation}
\label{eq:divFreeVel}
\begin{pmatrix} u\\ v \end{pmatrix} = \begin{pmatrix} \sum_{n,m} \frac{a_{nm}}{\lambda_{nm}}\phi_n \phi_m^y \\ \sum_{n,m} \frac{a_{nm}}{\lambda_{nm}}\phi_n^x \phi_m \end{pmatrix},
\end{equation}
where $\phi_n^x$ and $\phi_n^y$ are the derivatives of the complex exponential basis functions with respect to $x$ and $y$ respectively. This field is divergence-free.

Note that a periodic basis is used to generate the divergence-free field. As the image sequence we use depicts rotation within the domain with little movement on the boundary, this basis is suitable.  

\begin{figure}[htbp]
  \centering
  \includegraphics[width=0.6\textwidth]{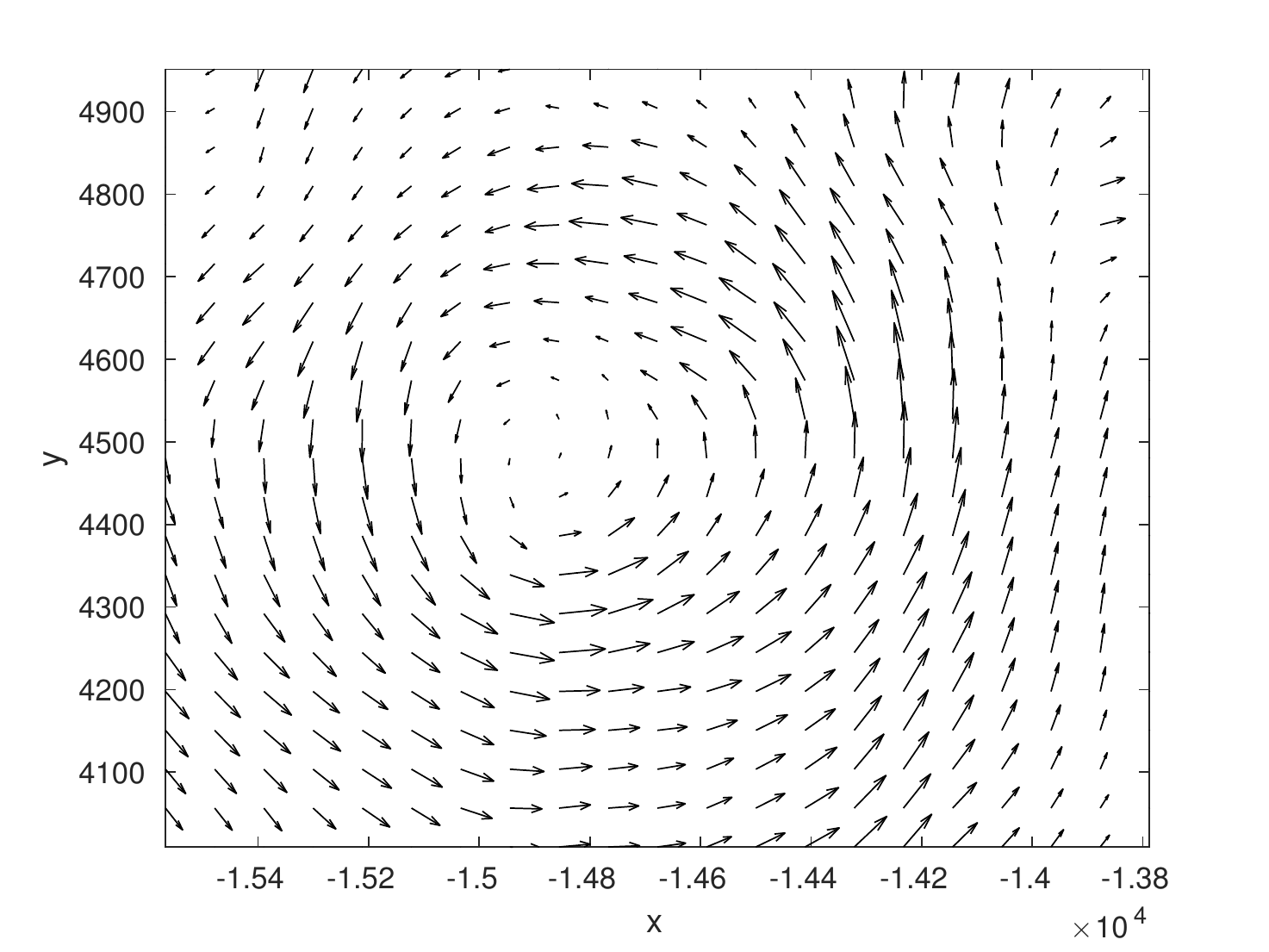}
  \caption{Sample velocity field for advection of satellite images}
  \label{fig:windField}
\end{figure}

The velocity field that advects the images is obtained using an optical flow estimation procedure on a similar set of images to the ones being assimilated. The reason we do not use the observations to obtain the optical flow velocity field is that we do not wish the state trajectory to precisely pass through the images; rather, the velocity field only roughly captures flow and is thus a source of uncertainty. This way, the images are used to steer the state of the model, which will veer off track due to the imprecise velocity field. The images are interpolated onto the DG-LGL grid using bi-linear interpolation. An example of the advection field is shown in~\cref{fig:windField}.

\subsubsection{Spatially and temporally sparse satellite observations}
\label{sec:CODchequered}

\begin{figure}
        \centering
        \begin{subfigure}[b]{0.49\textwidth}
                \includegraphics[width=\textwidth]{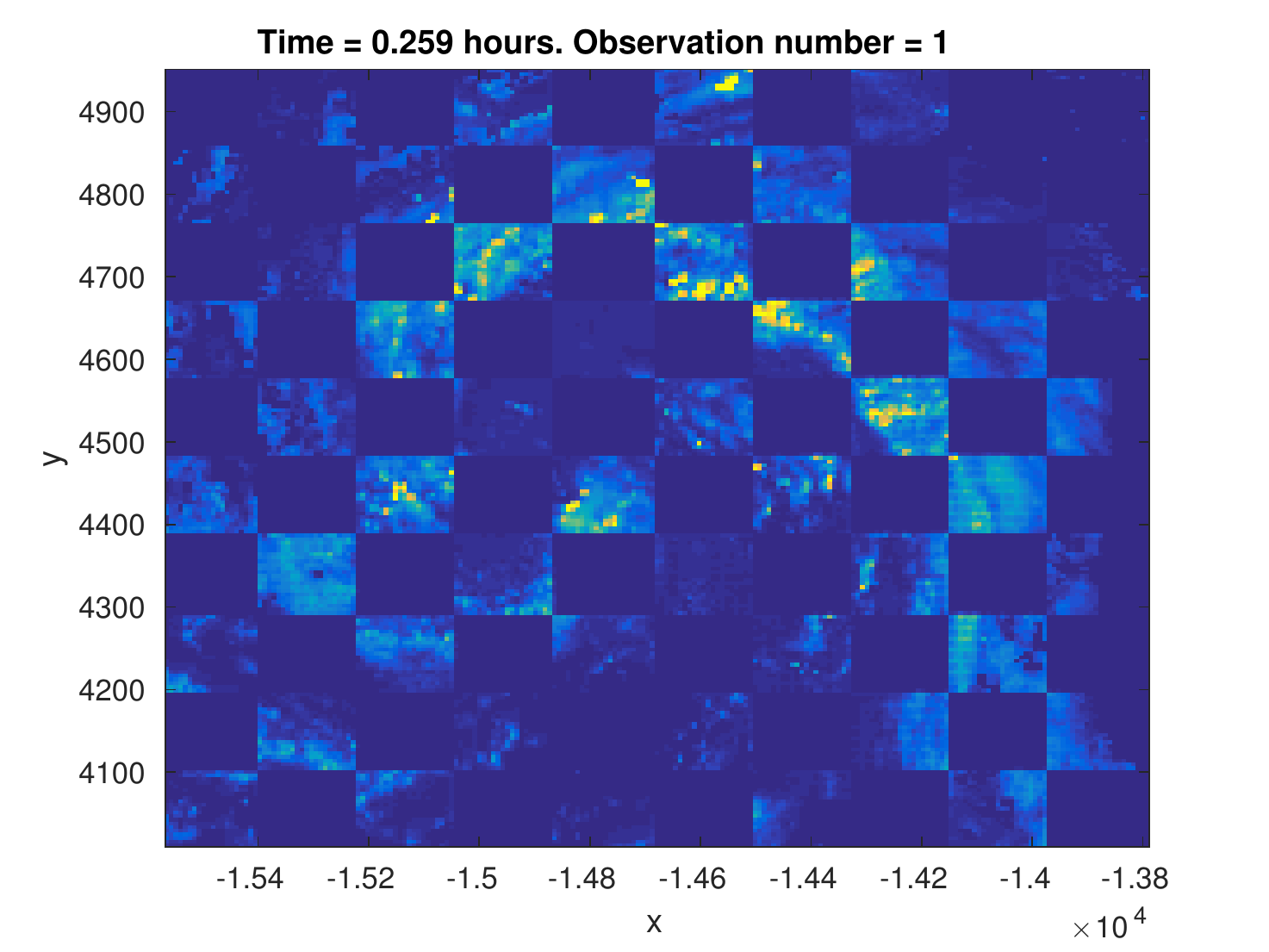}
                \caption{Assimilation of 1st observation (Relative error $\approx 0.56$)}
                \label{fig:localFilterPartialObsCOD1}
        \end{subfigure}%
        ~ 
        \begin{subfigure}[b]{0.49\textwidth}
                \includegraphics[width=\textwidth]{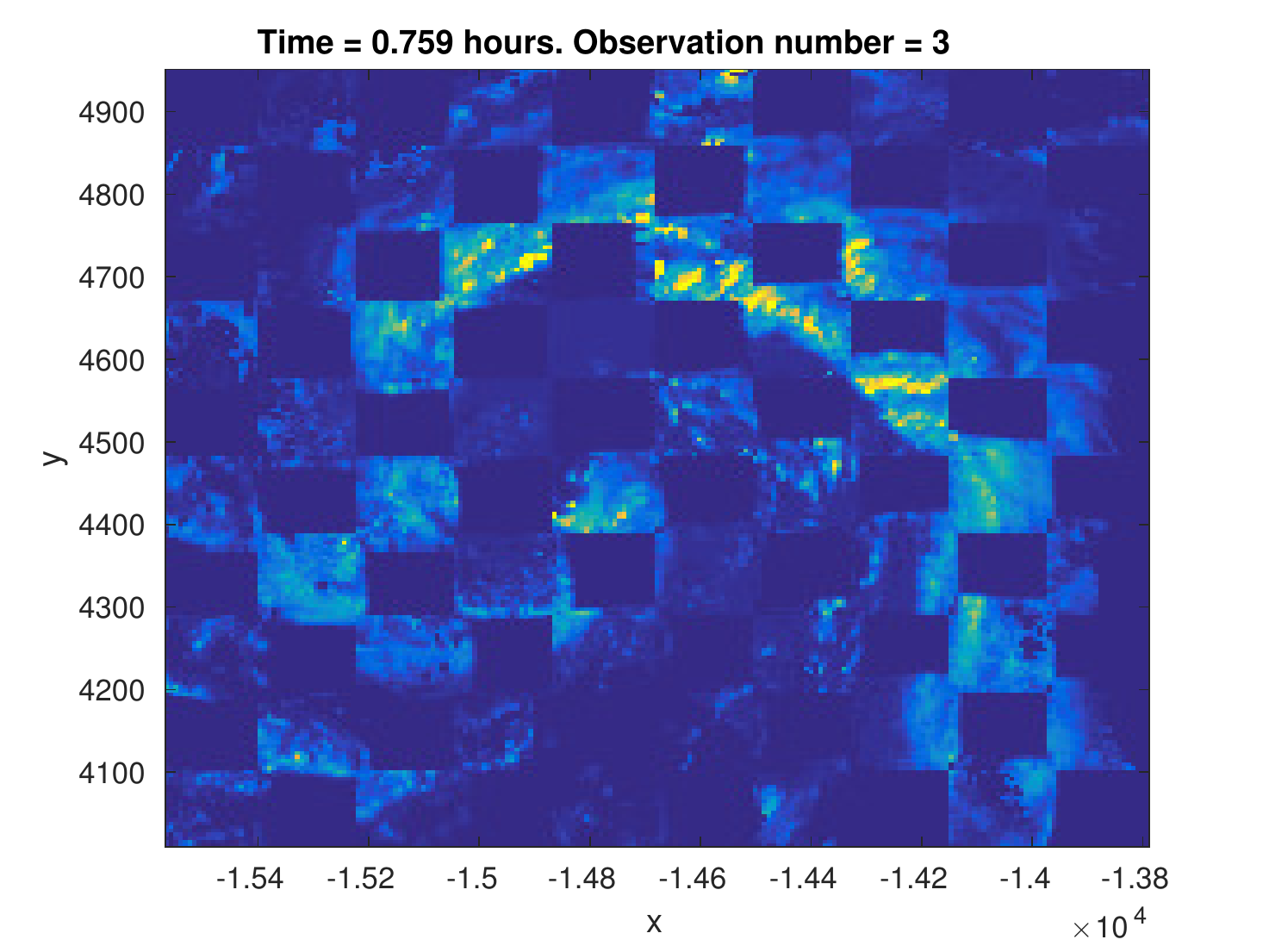}
                \caption{Assimilation of 3rd observation (Relative error $\approx 0.46$)}
                \label{fig:localFilterPartialObsCOD2}
        \end{subfigure}%

        \begin{subfigure}[b]{0.49\textwidth}
                \includegraphics[width=\textwidth]{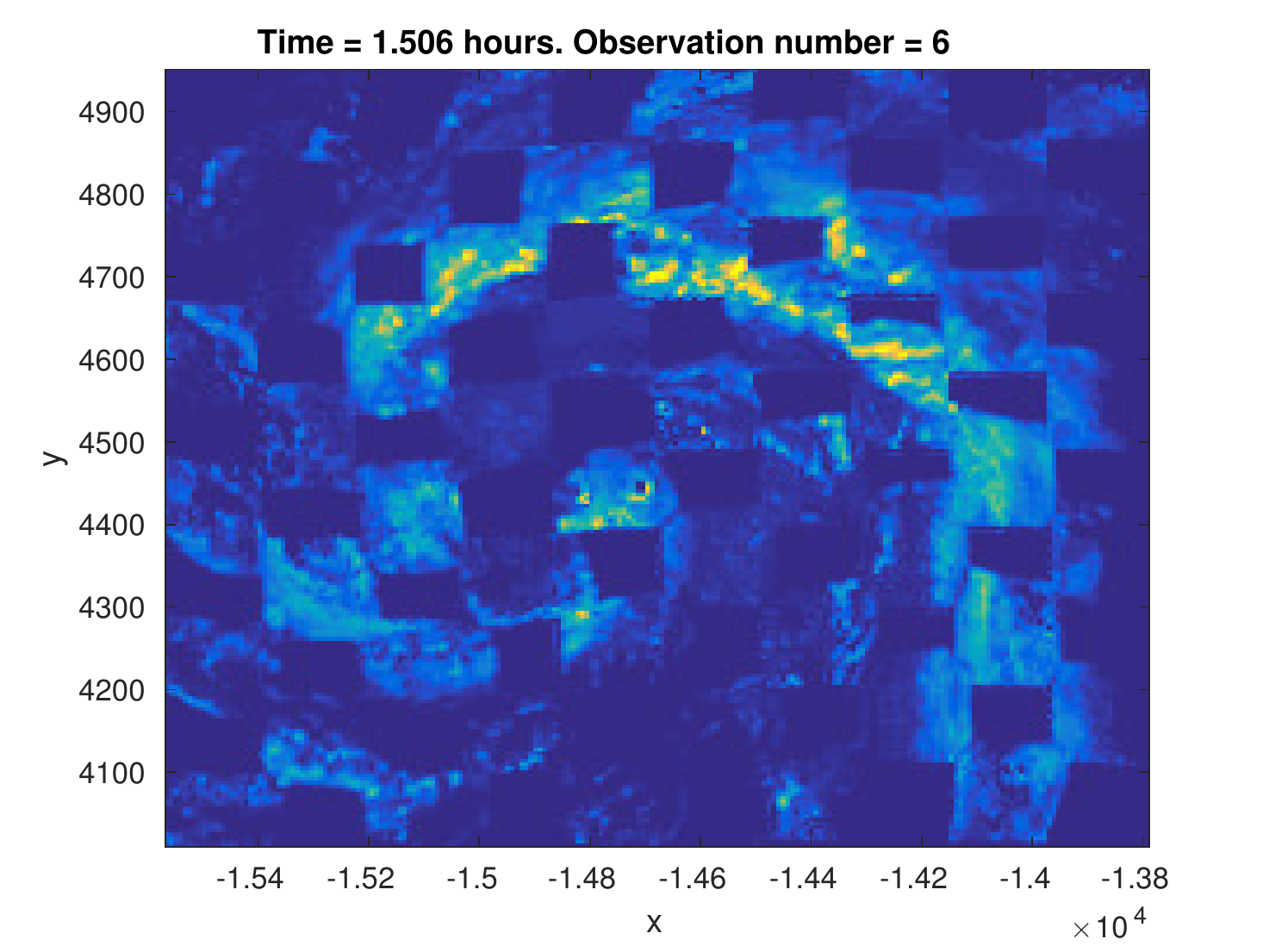}
                \caption{Assimilation of 6th observation (Relative error $\approx 0.35$)}
                \label{fig:localFilterPartialObsCOD3}
        \end{subfigure}%
        ~ 
        \begin{subfigure}[b]{0.49\textwidth}
                \includegraphics[width=\textwidth]{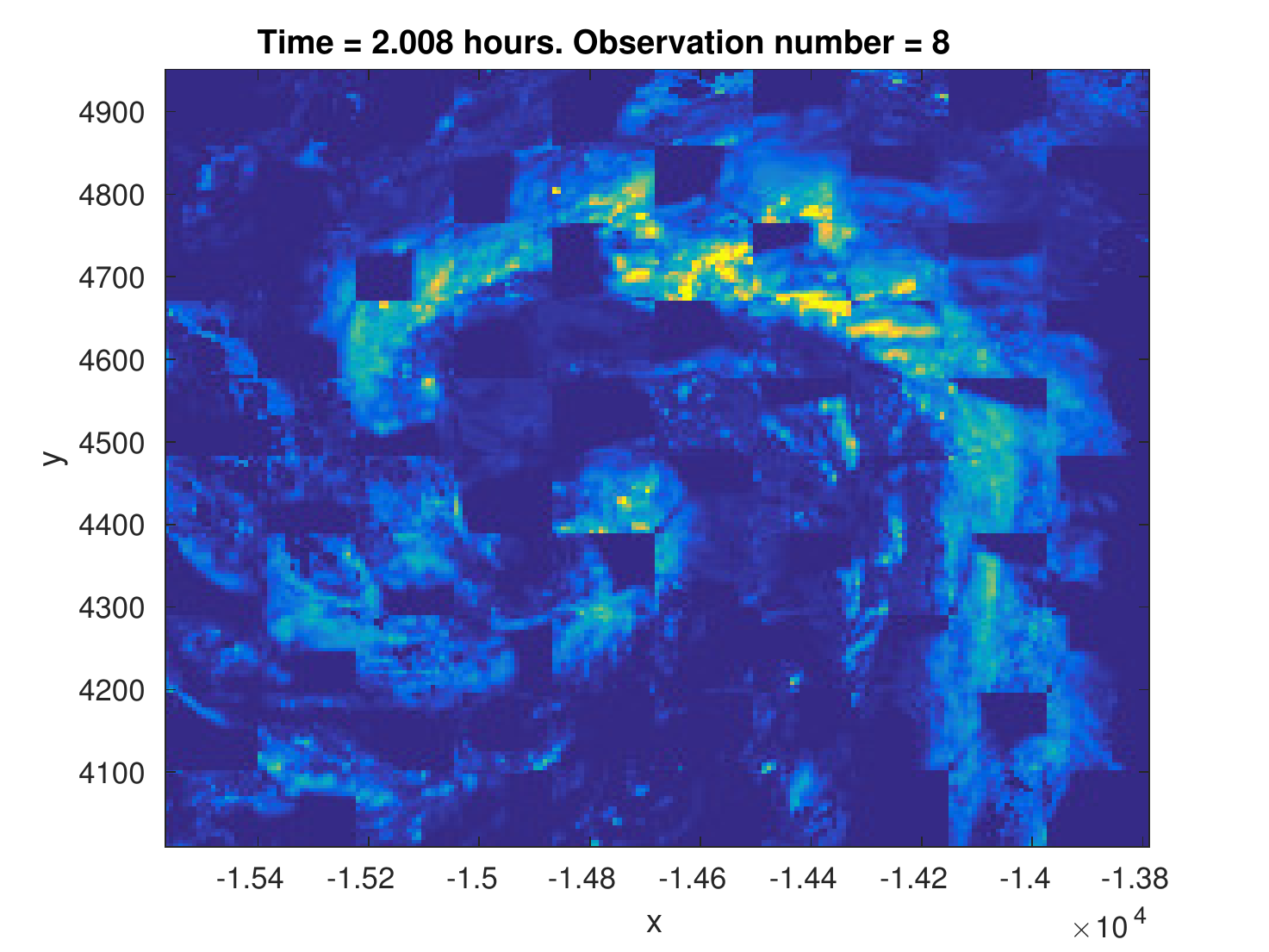}
                \caption{Assimilation of 8th observation (Relative error $\approx 0.28$)}
                \label{fig:localFilterPartialObsCOD4}
        \end{subfigure}%

        \begin{subfigure}[b]{0.49\textwidth}
                \includegraphics[width=\textwidth]{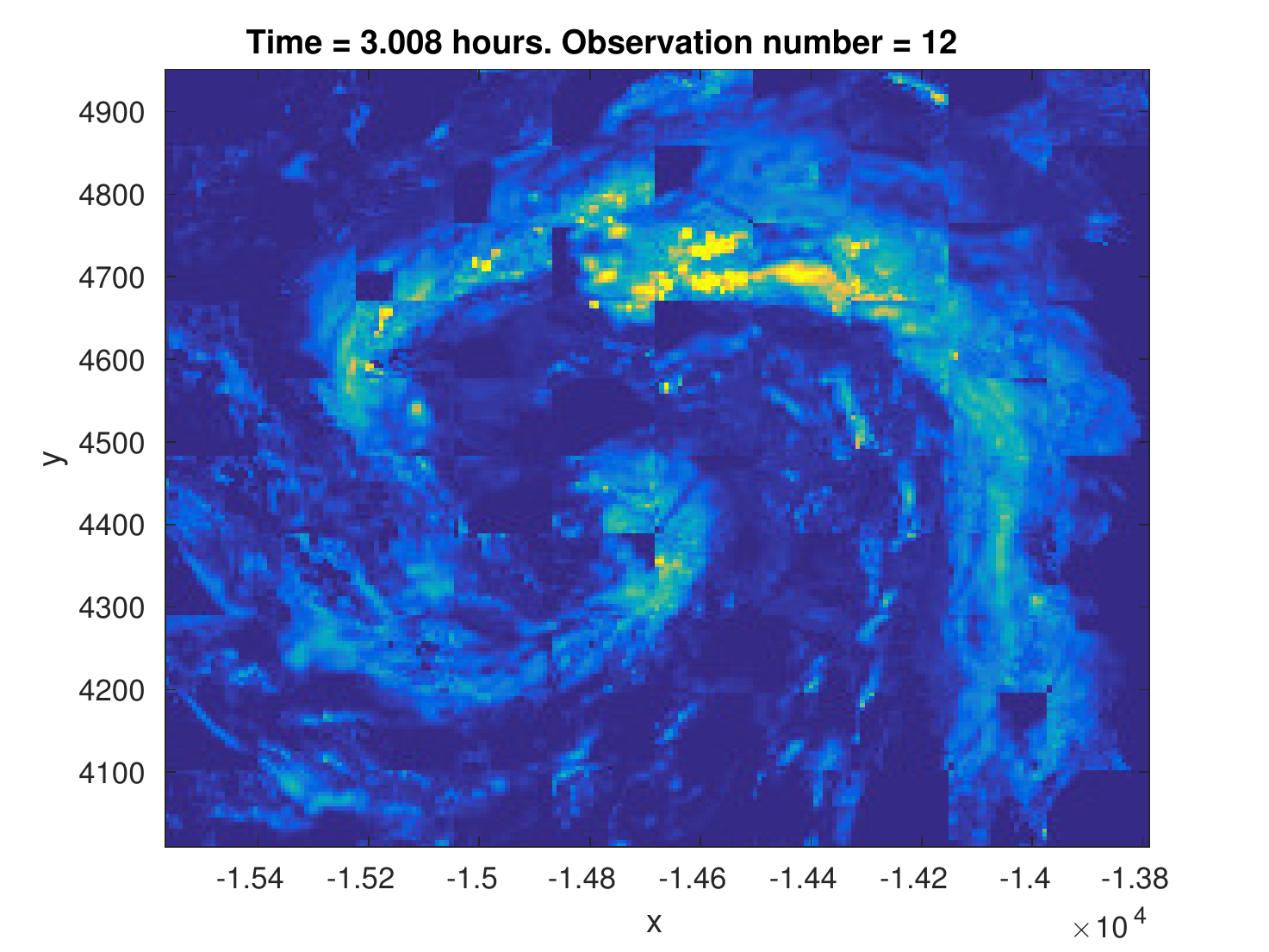}
                \caption{Assimilation of 12th observation (Relative error $\approx 0.20$)}
                \label{fig:localFilterPartialObsCOD5}
        \end{subfigure}%
        ~ 
        \begin{subfigure}[b]{0.49\textwidth}
                \includegraphics[width=\textwidth]{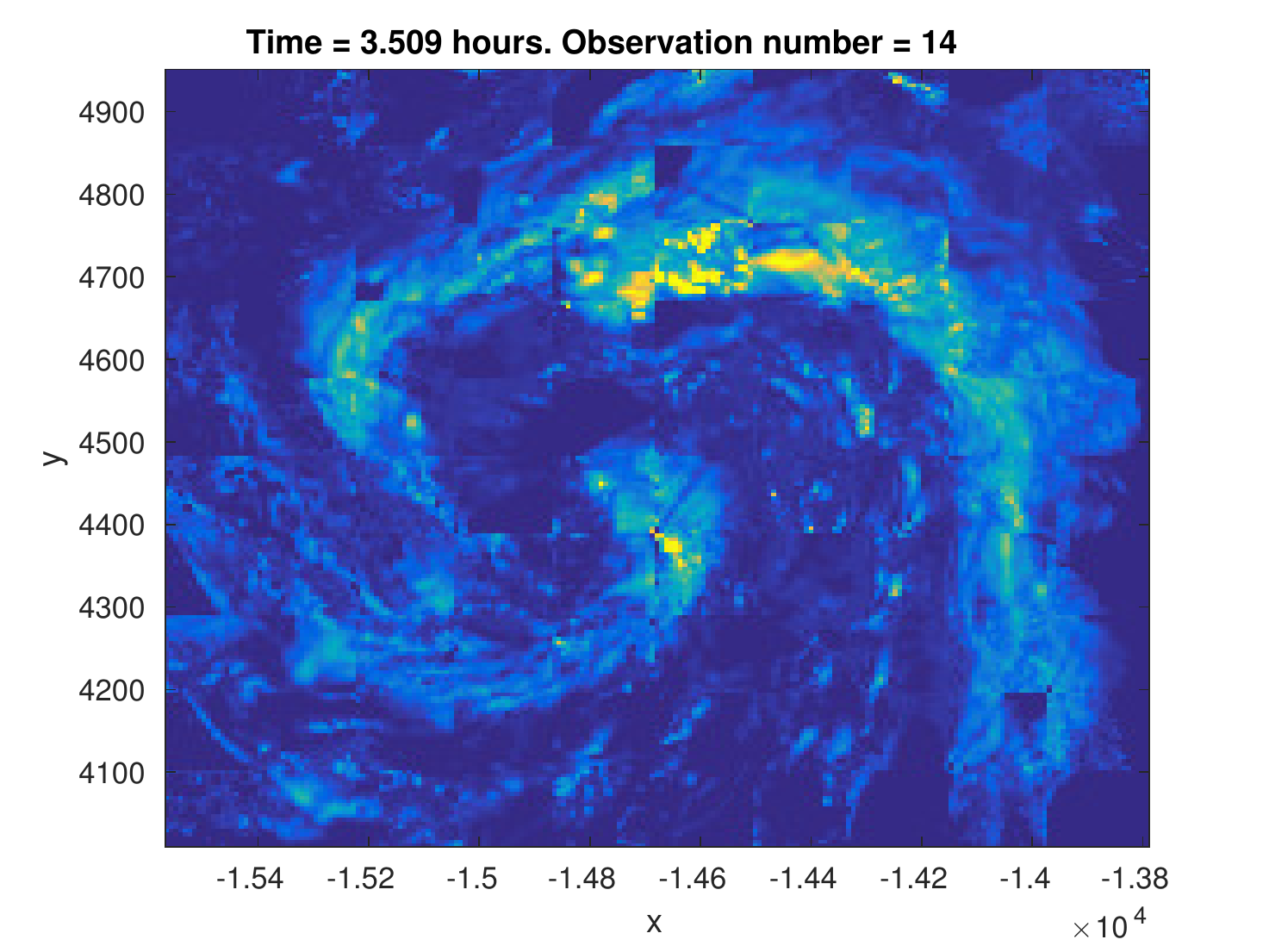}
                \caption{Assimilation of 14th observation (Relative error $\approx 0.16$)}
                \label{fig:localFilterPartialObsCOD6}
        \end{subfigure}%

    \caption{Assimilation of a partial satellite images with approximate divergence-free optical flow velocity field}\label{fig:localFilterPartialObsCOD}
\end{figure}

We perform an experiment to test the distributed filter with sparse observations by using the same ``chequered'' observation pattern that was used in~\cref{sec:sparseObsSynth} and \cref{sec:sparseObsSynthWrongBC} for the experiments on synthetic data in order to emulate partial image availability. However, instead of observing every other element, we observe the images on fixed $7 \times 7$ blocks of elements, separated by blocks of the same size, giving a $10\times 10$ chequered observation pattern. The time-step, $\Delta T$, computed for the $70 \times 70$ element grid with $N=3$, is adjusted so that it fits the 15-minute interval, yielding $\Delta T = 15/18$ minutes. Using this time-step, observations become available to the filter at intervals of 18 time-levels. This temporal sparsity can not be reduced for the current grid, as increasing $\Delta T$ so that the observations are temporally less sparse would cause the system to violate the CFL condition. For the experiment, we initialise the state to zero, i.e., $\bm{c}_0=0$, and impose zero Dirichlet boundary conditions at inflow nodes and free-exit Neumann conditions at outflow nodes. $q$, $q_0$ and $q_\partial$ are intialised as in~\cref{sec:sparseObsSynth} and $g=g_0=g_\partial=1$. The distributed filter is run for a simulated time of $3.5+$ hours over which time it assimilates 14 observations. In~\cref{fig:localFilterPartialObsCOD}, the results are shown for 6 of those observations including the first and the last.
The observation pattern is clear from~\cref{fig:localFilterPartialObsCOD1} where we see the alternating groups of $7 \times 7$ elements with/without available observations. Over time, the error decreases, with the last observation yielding a relative error of $\approx 0.16$. Over time, the observation pattern does not change, so the portions of the domain without data rely on flux from neighbouring elements as we also saw in synthetic experiments in~\cref{sec:sparseObsSynth} and \cref{sec:sparseObsSynthWrongBC}. We refer the reader to~\cite{SZTTAA_CDC17} for further details on this.

\section{Conclusions}
\label{sec:conclusions}

We conclude by summarising features of the proposed distributed filtering approach and outlining directions of the future research. The key advantage of the proposed filtering framework is its scalability: the computational cost of the distributed state estimator equals the cost of the filtering at each element multiplied by the total number of elements. Another important feature is the structure preserving discretisation at each element which preserves symmetry and positivity of the gain $P^k$, and all quadratic invariants of the estimation error dynamics. Finally, the filter is designed for the hyperbolic equation directly, without introducing artificial viscosity. The latter is, to some extent, introduced by means of the Lax-Friedrichs flux which ``smooths'' out the jumps across inter-elemental boundaries. Experimental assessment of the distributed filters on synthetic and real data shows great potential of the algorithm.

A key direction for the future work is to study convergence of the proposed scheme, and introduce explicit boundary interconnection mechanisms for the elemental Riccati equations to improve the uncertainty exchange mechanism.

\appendix
\section{Proofs}
\label{s:proofs}
\begin{proof}[Proof of \cref{p:1}]
Let $\bm{c}_h^k $ solve~\cref{eq:affineSys}. Then $\bm{c}_h^k = \bar{\bm{c}}_h^k+ \tilde{\bm c}_h^k$, provided $\bar{\bm{c}}_h^k$ solves $\dfrac{d \bar{\bm c}_h^k}{dt} = A^k(t) \bar{\bm c}_h^k + \bm{b}^k$, $ \bar{\bm c}_h^k(0)=\bm{c}^k_0$, and $\dfrac{d \tilde{\bm c}_h^k}{dt} = A^k(t) \tilde{\bm c}_h^k + V^k\tilde{\bm e}^k(t)$, $\tilde{\bm c}_h^k(0) = G_0^k\bm{e}^k_0$. Define $\tilde{\bm{y}}^k:=\bm{y}^k-\bar{\bm{c}}_h^k$, set $p:=(N+1)^2$, and let us find  $\hat {\bm{u}}\in L^2(t_0,t_f,\mathbb{R}^p)$ such that, $\forall \bm{u}\in L^2(t_0,t_f,\mathbb{R}^p)$: \[
\max_{\tilde{\bm{e}}^k,\bm{e}^k_0,E \bm{\eta}^k (\bm{\eta}^k)^\top }\sigma(\hat{\bm{u}})\le \max_{\tilde{\bm{e}}^k,\bm{e}^k_0,E \bm{\eta}^k (\bm{\eta}^k)^\top }\sigma(\bm{u})\,, \quad \sigma(\bm{u}):= E(\bm{\ell}^\top \tilde{\bm c}_h^k(t+s) - \int_{t}^{t+s} \bm{u}^\top\tilde{\bm{y}}^k d\tau)^2 \,,
\]
provided that~\cref{eq:ellips_Dk} and~\cref{eq:ellips_Y} hold true. Here $\bm{\ell}\in \mathbb{R}^p$ is some vector. Introducing  adjoint variable $\dfrac{d \bm{z}}{dt} = -(A^k)^\top \bm{z} + (H^k)^\top \bm{u}$, $\bm{z}(t+s) = \bm{\ell}$, and integrating by parts we find: \[
\sigma(\bm{u})=\left(\bm{z}^\top(t)(G_0^k\bm{e}^k_0) + \int_{t}^{t+s} \bm{z}^\top(V^k\tilde{\bm e}^k)d\tau\right)^2 + E\left(\int_t^{t+s} \bm{u}^\top\bm{\eta}^k d\tau\right)^2
\]
By using Cauchy-Schwartz inequality, \cref{eq:ellips_Dk} and~\cref{eq:ellips_Y} we find: \[
\max_{\tilde{\bm{e}}^k,\bm{e}^k_0,E \bm{\eta}^k (\bm{\eta}^k)^\top }\sigma(\bm{u}) = (1+2s)(N+1)^2 (\bm{z}^\top(t) G_0^k(Q_0^{k})^{-1}G_0^k \bm{z}(t) + \int_t^{t+1} \bm{z}^\top(V^k(\tilde{Q}^k)^{-1}(V^k)^\top)\bm{z} + \bm{u}^\top (R^k) \bm{u} d\tau)
\] By using standard LQ control theory results~\cite{ZhukCDC14} we find that the unique minimum point of this quadratic cost functional along the solutions of the adjoint equation for $\bm{z}$ satisfies the following feed-back representation: $\hat{\bm{u}} = (R^k)^{-1}H^kP^k\bm{z}$. By using the latter, and integration by parts it is not hard to find that: $\int_t^{t+s} \hat{\bm{u}}^\top \tilde{\bm{y}}^kd\tau = \bm{\ell}^\top \hat{\tilde{\bm c}}_h^k(t+s)$, provided $\hat{\tilde{\bm c}}_h^k$ solves \[
\frac{d \hat{\tilde{\bm{c}}}^k_h}{dt} = A^k\hat{\tilde{\bm{c}}}_h^{k} + P^k (H^k)^\top (R^k)^{-1}(\tilde{\bm{y}}^k - H^k\hat{\tilde{\bm{ c}}}^k_h)\,, \quad \hat{\tilde{\bm{ c}}}^k_h(0) = 0\,.
 \] Integrating by parts we find:
 \begin{equation}
   \label{eq:me}
\max_{\tilde{\bm{e}}^k,\bm{e}^k_0,E \bm{\eta}^k (\bm{\eta}^k)^\top }\sigma(\hat{\bm{u}})=  E(\bm{\ell}^\top \tilde{\bm c}_h^k(t+s) -  \bm{\ell}^\top \hat{\tilde{\bm c}}_h^k(t+s))^2 = (1+2s)(N+1)^2\bm{\ell}^\top P^k(t+s)\bm{\ell}
\end{equation}
Now, by recalling that $\bm{c}_h^k = \bar{\bm{c}}_h^k+ \tilde{\bm c}_h^k$, and by noting that, in fact:  $\hat{\bm{c}}^k_h = \bar{\bm{c}}_h^k+\hat{\tilde{\bm{c}}}^k_h$ we obtain~\cref{eq:lMF_err} from~\cref{eq:me} by setting $\bm{\ell}=(0,\dots, 0,1,0,\dots, 0)^\top$, where $1$ is at position $j$. This completes the proof.
\end{proof}

\bibliographystyle{siamplain}

\end{document}